\documentclass{amsart}
\usepackage{amsfonts}
\usepackage{amsmath}
\usepackage{amssymb}
\usepackage{hyphenat}
\usepackage{diagrams}
\usepackage[dvips]{graphicx}
\usepackage{float}

\newarrow {Mapsto} |--->
\newarrow {TeXto} ----{->}
\newarrow {TeXonto} ----{->>}
\newarrow {Line} -----
\diagramstyle[PostScript=dvips]

\newtheorem{theorem}{Theorem}[section]
\theoremstyle{plain}

\newtheorem{corollary}[theorem]{Corollary}

\newtheorem{lemma}[theorem]{Lemma}

\newtheorem{problem}[theorem]{Problem}

\newenvironment{proofnosquare}[1][Proof]{\smallskip\noindent\textsl{#1.} }{}
\newenvironment{Theorem Ref}[1][Theorem]{\medskip\noindent\textbf {#1}}{\medskip}

\makeatletter \@addtoreset{figure}{section} \makeatother

\newcommand{\arf}{\operatorname{Arf}}
\newcommand{\aut}{\operatorname{Aut}}
\newcommand{\coker}{\operatorname{coker}}
\newcommand{\Hom}{\operatorname{Hom}}
\newcommand{\im}{\operatorname{Im}}
\newcommand{\image}{\operatorname{Image}}
\newcommand{\lk}{\operatorname{lk}}
\newcommand{\Sp}{\operatorname{Sp}}
\newcommand{\tor}{\operatorname{Torsion}}

\begin{document}
\title{Bordism Invariants of the Mapping Class Group}
\author{Aaron Heap}
\date{November 10, 2004}

\begin{abstract}We define new bordism and spin bordism invariants of
certain subgroups of the mapping class group of a surface. In
particular, they are invariants of the Johnson filtration of the
mapping class group. The second and third terms of this filtration
are the well-known Torelli group and Johnson subgroup,
respectively. We introduce a new representation in terms of spin
bordism, and we prove that this single representation contains all
of the information given by the Johnson homomorphism, the
Birman-Craggs homomorphism, and the Morita homomorphism.
\end{abstract}

\maketitle

\section{Introduction}

Let $\Sigma _{g,1}$ be a compact, oriented surface of genus $g$
with one boundary component. Let $\Gamma _{g,1}$ be the mapping
class group of $\Sigma _{g,1}$. That is, $\Gamma _{g,1}$ is the
group of isotopy classes of orientation-preserving homeomorphisms
of $\Sigma _{g,1}$ which fix the boundary. The study of mapping
class groups has important applications in many different areas of
topology, differential geometry, and algebraic geometry. Here we
are particularly interested in $\Gamma _{g,1}$ within the area of
3-manifold topology.

The mapping class group $\Gamma _{g,1}$ acts naturally by
automorphisms on the fundamental group $F=\pi _{1}(\Sigma
_{g,1})$, which is a free group of rank $2g$. Then we have the
induced representation $\Gamma_{g,1}\rightarrow \aut(F)$, and this
representation is known classically to be injective. Let $\left\{
F_{k}\right\} _{k\geq 1}$ be the lower central series of $F$. That
is, $F_{1}=F$ and the rest of the terms are defined inductively by
$F_{k+1}=\left[ F_{k},F_{1}\right] $ for any $k\geq 1$. Then
$\Gamma _{g,1}$ acts naturally on the nilpotent quotients
$F/F_{k}$, providing a series of representations
\begin{equation*}
\rho _{k}:\Gamma _{g,1}\rightarrow \aut \left(
\frac{F}{F_{k}}\right) .
\end{equation*}%
Note that $F/F_{2}$ is isomorphic to the first homology group
$H_{1}=H_{1}(\Sigma _{g,1};\mathbb{Z)}$, and $\rho _{2}$ is the
same as the classical representation $\Gamma _{g,1}\rightarrow
\Sp(2g;\mathbb{Z)}$ of the mapping class group onto the Siegel
modular group, which is the group of symplectic automorphisms of
$H_{1}$ with respect to the skew-symmetric intersection pairing.

The \textit{generalized Johnson subgroup} $\mathcal{J}(k)\subseteq
\Gamma _{g,1}$ is defined to be the kernel of $\rho _{k}$. That
is, $\mathcal{J}(k)$ is the subgroup of the mapping class group
consisting of those homeomorphisms which induce the identity on
$F/F_{k}$. The subgroup $\mathcal{J}(2)=\mathcal{T}_{g,1}$ is more
commonly known as the\textit{\ Torelli group}, and
$\mathcal{J}(3)=\mathcal{K}_{g,1}$ is traditionally referred to as
the \textit{Johnson subgroup}. The Johnson subgroup was originally
defined to be the subgroup of $\Gamma _{g,1}$ generated by all
Dehn twists about separating simple closed curves on $\Sigma
_{g,1}$. The fact that these two definitions of
$\mathcal{K}_{g,1}$ are equivalent was proved by D. Johnson in
\cite{[J5]}.

To get a better understanding of the structure of the subgroup
$\mathcal{J}(k)$, it is natural to seek abelian representations
for it. That is, we would hope to understand $\mathcal{J}(k)$
better by investigating abelian quotients of it. The first such
quotient of the Torelli group $\mathcal{J}(2) $ was given by a
homomorphism due to D. Sullivan in \cite{[Su]}. Johnson gave
another homomorphism for $\mathcal{J}(2)$, of which Sullivan's is
a quotient, in \cite{[J2]}. He later generalized this homomorphism
to $\mathcal{J}(k)$ for all $k\geq 2$ in \cite{[J3]}, thus giving
a family of homomorphisms
\begin{equation*}
\tau _{k}:\mathcal{J}(k)\rightarrow \Hom \left(
H_{1},\frac{F_{k}}{F_{k+1}}\right) ,
\end{equation*}%
now known as the \textit{Johnson homomorphisms}. In the case
$k=2$, the image of $\tau _{2}$ is known to be a submodule
$D_{2}(H_{1})$ of $\Hom(H_{1},F/F_{2})=\Hom(H_{1},H_{1}).$
Moreover, the kernel of $\tau _{2}$ is known to be
$\mathcal{J}(3)$. In general, $\ker \tau _{k}=\mathcal{J}\left(
k+1\right)$. However, the image of $\tau _{k}$ is not known for
$k\geq 3$, and it is a fundamental problem in the study of the
mapping class group to determine its image.

In \cite{[BC]} J. Birman and R. Craggs produced a collection of
abelian quotients of $\mathcal{J}(2)$ given by homomorphisms onto
$\mathbb{Z}_{2}, \rho :\mathcal{J}(2)\rightarrow \mathbb{Z}_{2}.$
These are finite in number and unrelated to Johnson's
homomorphism. However, Johnson showed in \cite{[J6]} that the
Johnson homomorphism $\tau _{2}$ and the totality of these
\textit{Birman-Craggs homomorphisms}, together, completely
determine the abelianization of the Torelli group $\mathcal{J}(2)$
for $g\geq 3$. The abelianization of $\mathcal{J}(k)$ is not known
for $k>2$.

In this paper we give new representations in terms of the
3-dimensional bordism groups $\Omega _{3}(F/F_{k})$ and $\Omega
_{3}^{spin}(F/F_{k})$. The former is a faithful representation of
the abelian quotient
$\mathcal{J}(k)/\mathcal{J}\left(2k-1\right)$, and the latter is a
homomorphism which combines the Johnson and Birman-Craggs
homomorphisms into a single homomorphism.  See Sections
\ref{(BRMCG)} and \ref{(SBRMCG)} for specific details.

\section{The Johnson Homomorphism\label{(JH)}}

\subsection{Johnson's Original Definition of $\protect\tau
_{k}$\label{(JOD)}}

In this section we give a description of Johnson's homomorphisms.
Let $\Sigma _{g,1}$ be a compact, oriented surface of genus $g$
with one boundary component and with fundamental group $F$. Let
$\left\{ F_{k}\right\} _{k\geq 1}$ be the lower central series of
$F$. Let the \textit{generalized Johnson subgroup}
$\mathcal{J}(k)$ be the subgroup of the mapping class group
consisting of those homeomorphisms that induce the identity on
$F/F_{k}$.

Consider any $f\in \mathcal{J}(k)$. Choose a representative
$\gamma \in \pi _{1}(\Sigma _{g,1})=F$ for any given element
$\left[ \gamma \right] \in H_{1}=H_{1}(\Sigma
_{g,1};\mathbb{Z)}=F/F_{2}$, and consider the element $f_{\ast
}(\gamma )\gamma ^{-1}$ which belongs to $F_{k}$ since $f\in
\mathcal{J}(k)$ implies $f_{\ast }$ acts trivially on $F/F_{k}$.
Then let $\left[ f_{\ast }(\gamma )\gamma ^{-1}\right] \in
F_{k}/F_{k+1}$ denote the equivalence class of $f_{\ast }(\gamma
)\gamma ^{-1}$ under the projection $F_{k}\rightarrow
F_{k}/F_{k+1}$. Then we define the \textit{Johnson homomorphisms}
\begin{equation*}
\tau _{k}:\mathcal{J}(k)\rightarrow \Hom \left(
H_{1},\frac{F_{k}}{F_{k+1}}\right)
\end{equation*}%
by letting $\tau _{k}(f)$ be the homomorphism $\left[ \gamma
\right] \rightarrow \left[ f_{\ast }(\gamma )\gamma ^{-1}\right]
$. The skew-symmetric intersection pairing on $H_{1}$ defines a
canonical isomorphism $H_{1}\cong \Hom(H_{1},\mathbb{Z})$, and
this induces an isomorphism
\begin{equation*}
\Hom \left( H_{1},\frac{F_{k}}{F_{k+1}}\right) \cong
\Hom(H_{1},\mathbb{Z})\otimes \frac{F_{k}}{F_{k+1}}\cong
H_{1}\otimes \frac{F_{k}}{F_{k+1}}.
\end{equation*}%
Thus we could also write
\begin{equation*}
\tau _{k}:\mathcal{J}(k)\rightarrow H_{1}\otimes
\frac{F_{k}}{F_{k+1}}\text{.}
\end{equation*}

This is Johnson's original definition \cite{[J3]}, but there are
several equivalent definitions of his homomorphism. Also in
\cite{[J3]}, one can see a definition in terms of the intersection
ring of the mapping torus of $f$. There is a definition of $\tau
_{k}$ in terms of the Magnus representation of the mapping class
group $\Gamma _{g,1}$ that may be found in \cite{[Ki]} or
\cite{[Mo]}.

The final definition we mention in this paper will be given in
Section \ref{(MPD)}. It was stated by Johnson \cite{[J3]} and
verified by T. Kitano \cite{[Ki]}. This definition gives a
computable description of $\tau _{k}$ in terms of Massey products
of mapping tori.

We complete this section with a few well-known facts about the
Johnson homomorphisms $\tau _{k}$ and the subgroups
$\mathcal{J}(k)$. It was shown by Morita in \cite{[Mo]} that
\begin{equation*}
\left[ \mathcal{J}(k),\mathcal{J(}l)\right] \subset \mathcal{J(}k+l-1).
\end{equation*}%
In particular, the commutator subgroup $\left[
\mathcal{J}(k),\mathcal{J}(k)\right]$ is a subgroup of
$\mathcal{J(}2k-1)$ for $k\geq 2$. As mentioned before, $\ker \tau
_{k}=\mathcal{J}\left(k+1\right)$. Then the image of $\tau _{k}$
is isomorphic to the abelian quotient
$\mathcal{J}(k)/\mathcal{J}\left(k+1\right)$. Thus the information
provided by the $k-1$ homomorphisms $\tau _{k},...,\tau _{2k-2}$
can be combined to determine the abelian quotient
$\mathcal{J}(k)/\mathcal{J}\left(2k-1\right).$ Unfortunately this
only at most detects the free-abelian part of the abelianization
$\mathcal{J}(k)/\left[ \mathcal{J}(k),\mathcal{J}(k)\right] \cong
H_{1}(\mathcal{J}(k))$. For example, the image of $\tau _{2}$ is
given by
$\mathcal{J}(2)/\mathcal{J}(3)=\mathcal{T}_{g,1}/\mathcal{K}_{g,1},$
and $\mathcal{J}(2)/\mathcal{J}(3)\otimes \mathbb{Q}\cong
H_{1}(\mathcal{T}_{g,1};\mathbb{Q)}$, whereas the abelianization
of the Torelli group $H_{1}(\mathcal{T}_{g,1})$ has 2-torsion. We
will discuss this 2-torsion in more detail in Section \ref{(BCH)}.

\subsection{Massey Products\label{(MP)}}

Let $\left( X,A\right) $ be a pair of topological spaces, and
unless otherwise stated we assume that the coefficients for
homology and cohomology groups are always the integers
$\mathbb{Z}$. In this section we will give the definition of the
Massey product
\begin{equation*}
H^{1}(X,A)\otimes \cdot \cdot \cdot \otimes H^{1}(X,A)\rightarrow H^{2}(X,A)
\end{equation*}%
since these are the only dimensions that we are interested in
using, and we will give a few useful properties of which we wish
to take advantage. The general definition is completely analogous
except for various sign conventions, and we refer the reader to D.
Kraines \cite{[Kr]}. For a more complete description of this
specific definition we are giving and for some useful examples, we
refer you to R. Fenn's book \cite{[Fe]}.

Massey products may be viewed as higher order analogues of cup
products and are defined when certain cup products vanish. Let
$u_{1},...,u_{n}\in H^{1}(X,A)$ be cohomology classes with cocycle
representatives $a_{1},...,a_{n}\in C^{1}(X,A)$, respectively.
A\textit{\ defining set }for the Massey product $\left\langle
u_{1},...,u_{n}\right\rangle $ is a collection of cochains
$a=\left( a_{i,j}\right) $, $1\leq i\leq j\leq n$ and $\left(
i,j\right) \neq \left( 1,n\right) $, satisfying

\begin{enumerate}
\item[(1)] $a_{i,i}=a_{i}$ for any $i\in \left\{ 1,...,n\right\} $,

\item[(2)] $a_{i,j}\in C^{1}(X,A)$,

\item[(3)] $\delta a_{i,j}=\sum_{r=i}^{j-1}a_{i,r}\cup a_{r+1,j}$.
\end{enumerate}

For such a defining set $a$ consider the cocycle $u(a)\in
C^{2}(X,A)$ given by
\begin{equation*}
u(a)=\sum_{r=1}^{n-1}a_{1,r}\cup a_{r+1,n}\text{.}
\end{equation*}%
The \textit{Massey product} $\left\langle
u_{1},...,u_{n}\right\rangle $ is defined if a defining set $a$
exists, and it is defined to be the subset of $H^{2}(X,A)$
consisting of the values $u(a)$ of all such defining sets $a$.

The length 1 Massey product $\left\langle u_{1}\right\rangle $ is
simply defined to be $u_{1}$, and its defining set is any cocycle
representative of $u_{1}$. The length 2 Massey product
$\left\langle u_{1},u_{2}\right\rangle $ is the cup product
$u_{1}\cup u_{2}$. The triple Massey product $\left\langle
u_{1},u_{2},u_{3}\right\rangle $ is defined only when
$\left\langle u_{1},u_{2}\right\rangle $ and $\left\langle
u_{2},u_{3}\right\rangle $ are zero. As you may notice from the
definition, Massey products of length 3 or greater may not be
uniquely defined but in fact may be a set of elements. However, if
a sufficient number of smaller Massey products vanish, then
$\left\langle u_{1},...,u_{n}\right\rangle $ is uniquely defined.
We have the following useful properties.

(\ref{(MP)}.1) \textit{Uniqueness.} For $n\geq 3$, the Massey
product $\left\langle u_{1},...,u_{n}\right\rangle $ is uniquely
defined if all Massey products of length less than $n$ are defined
and vanish. (This hypothesis is stronger than necessary for
uniqueness, but it is sufficient for our purposes.)

(\ref{(MP)}.2) \textit{Naturality.} Let $\left( Y,B\right) $ be a
pair of topological spaces, and consider a map of pairs $f:\left(
Y,B\right) \rightarrow \left( X,A\right) $. If $\left\langle
u_{1},...,u_{n}\right\rangle $ is defined then so is $\left\langle
f^{\ast }(u_{1}),...,f^{\ast }(u_{n})\right\rangle $, and $f^{\ast
}\left\langle u_{1},...,u_{n}\right\rangle \subset \left\langle
f^{\ast }(u_{1}),...,f^{\ast }(u_{n})\right\rangle $. Furthermore,
if $f^{\ast }$ is an isomorphism, then equality holds.

\subsection{Massey Product Description of $\protect\tau
_{k}$\label{(MPD)}}

We are now prepared to describe Johnson's homomorphisms $\tau
_{k}$ using Massey products of mapping tori. For a more complete
description, see the work of Kitano \cite{[Ki]}. As before,
$\Sigma _{g\text{,}1}$ is an oriented surface of genus
$g$\textit{\ }with one boundary component $\partial \Sigma
_{g,1}$. Consider any homeomorphism $f\in \mathcal{J}(k)$, and let
$T_{f,1}$ denote the mapping torus of $f$. That is, $T_{f,1}$ is
$\Sigma _{g,1}\times \left[ 0,1\right] $ with $x\times \left\{
0\right\} $ glued to $f(x)\times \left\{ 1\right\} $. Note that
the boundary\ $\partial T_{f,1}$ is the torus $\partial \Sigma
_{g,1}\times S^{1}$. With the natural orientation on $\left[
0,1\right] $, we have a local orientation on $T_{f,1}$ given by
the product orientation. Moreover, since $f\in \mathcal{J}(k)$
acts trivially on $H_{1}=H_{1}(\Sigma _{g,1})$ as long as $k\geq
2$, the mapping torus $T_{f,1}$ is an oriented homology $\Sigma
_{g,1}\times S^{1}$, but the Massey product structure may be
different than that of $\Sigma _{g,1}\times S^{1}$.

First, fix a basis $\left\{ \alpha _{1},...,\alpha _{2g}\right\} $
for the free group $F=\pi _{1}(\Sigma _{g,1})$. Then if $\gamma $
represents a generator of $\pi _{1}(S^{1})$, we get the following
presentation of $\pi _{1}(T_{f,1})$:
\begin{equation*}
\pi _{1}(T_{f,1})=\left\langle \alpha _{1},...,\alpha _{2g},\gamma
\mid \left[ \alpha _{1},\gamma \right] f_{\ast }(\alpha
_{1})\alpha _{1}^{-1},..., \left[ \alpha _{2g},\gamma \right]
f_{\ast }(\alpha _{2g})\alpha _{2g}^{-1}\right\rangle .
\end{equation*}%
By denoting the homology classes of $\alpha _{i}$ and $\gamma $ by
$x_{i}$ and $y$, respectively, we obtain a basis for
$H_{1}(T_{f,1})$:
\begin{equation*}
\left\{ x_{1},...,x_{2g},y\right\} \in H_{1}(T_{f,1})\text{.}
\end{equation*}%
Then since $H^{1}(T_{f,1})\cong \Hom(H_{1}(T_{f,1}),\mathbb{Z})$,
we have a dual basis for $H^{1}(T_{f,1})$:
\begin{equation*}
\left\{ x_{1}^{\ast },...,x_{2g}^{\ast },y^{\ast }\right\} \in
H^{1}(T_{f,1}) \text{.}
\end{equation*}%
Let $j:\left( T_{f,1},\varnothing \right) \rightarrow \left(
T_{f,1},\partial T_{f,1}\right) $ be the inclusion map. The long
exact sequence of a pair shows $j_{\ast}:H_{1}(T_{f,1})\rightarrow
H_{1}(T_{f,1},\partial T_{f,1})$ has kernel generated by $y$. So
we have a basis for $H_{1}(T_{f,1},\partial T_{f,1})$:
\begin{equation*}
\left\{ j_{\ast }(x_{1}),...,j_{\ast }(x_{2g})\right\} \in
H_{1}(T_{f,1},\partial T_{f,1})\text{.}
\end{equation*}%
And this gives a corresponding basis for $H_{2}(T_{f,1})\cong
H^{1}(T_{f,1},\partial T_{f,1})$:
\begin{equation*}
\left\{ X_{1},...,X_{2g}\right\} \in H_{2}(T_{f,1})\text{.}
\end{equation*}

Let $\varepsilon :\mathbb{Z}\left[ F\right] \rightarrow
\mathbb{Z}$ be the augmentation map and let
\begin{equation*}
\frac{\partial }{\partial \alpha _{i}}:\mathbb{Z}\left[ F\right] \rightarrow
\mathbb{Z}\left[ F\right] \text{, }1\leq i\leq 2g
\end{equation*}%
be the Fox's free derivatives. Here $\mathbb{Z}\left[ F\right] $
is the integral group ring of the free group $F$. Finally, let
$\mathfrak{X}$ denote the ring of formal power series in the
noncommutative variables $t_{1},...,t_{2g}$, and let
$\mathfrak{X}_{k}$ denote the submodule of $\mathfrak{X}$
corresponding to the degree $k$ part. One can show $F_{k}/F_{k+1}$
is a submodule of $\mathfrak{X}_{k}$, where the inclusion map is
induced by
\begin{equation*}
F_{k}\ni \zeta \longmapsto \sum_{j_{1},...,j_{k}}\varepsilon \frac{\partial
}{\partial \alpha _{j_{1}}}\cdot \cdot \cdot \frac{\partial }{\partial
\alpha _{j_{k}}}(\zeta )t_{j_{1}}...t_{j_{k}}\in \mathfrak{X}_{k}\text{.}
\end{equation*}%
Then we have the following theorem.

\begin{theorem}[Kitano]
\label{(kitano)}There is a homomorphism
$\tau_{k}:\mathcal{J}(k)\rightarrow \Hom(H_{1},\mathfrak{X}_{k})$
defined by letting $\tau_{k}(f)$ be the homomorphism
\begin{equation*}
x_{i}\longmapsto \sum_{j_{1},...,j_{k}}\left\langle \left\langle
x_{j_{1}}^{\ast },...,x_{j_{k}}^{\ast }\right\rangle
,X_{i}\right\rangle t_{j_{1}}...t_{j_{k}}
\end{equation*}%
where $\left\langle \text{\ \ ,\ \ }\right\rangle $ is the dual
pairing of $H^{2}(T_{f,1})$ and $H_{2}(T_{f,1})$. Moreover, this
homomorphism is the same as the Johnson homomorphism.
\end{theorem}

The canonical restriction $H^{\ast }(T_{f,1},\partial
T_{f,1})\rightarrow H^{\ast }(T_{f,1})$ leads to the following
theorem that gives a relation between the algebraic structure of
the mapping class group $\Gamma _{g,1}$ and the topological
structure of the mapping torus $T_{f,1}$.

\begin{theorem}[Kitano]
\label{(kitcor)}For any $f\in \Gamma _{g,1}$, $f\in
\mathcal{J}\left(k+1\right)$ if and only if all Massey products of
length $m$ of
\begin{equation*}
H^{1}(T_{f,1},\partial T_{f,1})\otimes \cdot \cdot \cdot \otimes
H^{1}(T_{f,1},\partial T_{f,1})\rightarrow H^{2}(T_{f,1},\partial
T_{f,1})\rightarrow H^{2}(T_{f,1})
\end{equation*}%
vanish for any $m$ with $1<m\leq k$.
\end{theorem}

\subsection{Morita's Refinement of $\protect\tau _{k}$\label{(MR)}}

In this section we point out the work of Morita in \cite{[Mo]},
where Johnson's homomorphism $\tau _{k}$ was refined so as to
narrow the range of $\tau _{k}$ to a submodule $D_{k}(H_{1})$ of
$H_{1}\otimes F_{k}/F_{k+1}$. This enhancement is obtained via a
homomorphism
\begin{equation*}
\tilde{\tau}_{k}:\mathcal{J}(k)\rightarrow H_{3}\left( \frac{F}{F_{k}}\right)
\end{equation*}%
defined below. Recall that the homology of a group $G$ is
$H_{i}(G)\equiv H_{i}(K(G,1),\mathbb{Z)}$, where $K(G,1)$ is an
Eilenberg-MacLane space. (We determine the kernel of Morita's
refinement in Corollary \ref{(morita kernel)} below.)

Let $\zeta \in \pi _{1}(\Sigma _{g,1})=F$ represent\ the homotopy
class of a simple closed curve on $\Sigma _{g,1}$ parallel to the
boundary $\partial \Sigma _{g,1}$. Now we choose a 2-chain $\sigma
\in C_{2}(F)$ such that $\partial \sigma =-\zeta $. Since any
$f\in \Gamma _{g,1}$ is required by definition to fix the
boundary, we have $\partial (\sigma -f_{\#}(\sigma ))=-\zeta
-\left( -\zeta \right) =0.$ Thus $\sigma -f_{\#}(\sigma )$ is a
2-cycle. Because $H_{2}(F)$ is trivial, there is a 3-chain
$c_{f}\in C_{3}(F) $ such that $\partial c_{f}=\sigma
-f_{\#}(\sigma ).$ Note that, essentially, this is just a mapping
cylinder construction. Let $\bar{c}_{f}$ denote the image of
$c_{f}$ in $C_{3}(F/F_{k}).$ If $f\in \mathcal{J}(k)$ then
$f_{\#}$ acts as the identity on $F/F_{k}$. Thus we have $\partial
\bar{c}_{f}=\overline{\sigma -f_{\#}(\sigma
)}=\bar{\sigma}-f_{\#}(\bar{\sigma})=0 $, and $\bar{c}_{f}$ is a
3-cycle. Finally define $\left[ \bar{c}_{f} \right] \in
H_{3}(F/F_{k})$ to be the corresponding homology class, and we
define Morita's homomorphism
$\tilde{\tau}_{k}:\mathcal{J}(k)\rightarrow H_{3}(F/F_{k})$ to be
$\tilde{\tau}_{k}(f)=\left[ \bar{c}_{f}\right] $. It is shown in
\cite{[Mo]} that the homology class $\left[ \bar{c}_{f}\right] $
does not depend on the choices that were made, and we refer you
there for the details.

Now consider the extension
\begin{equation*}
0\rightarrow \frac{F_{k}}{F_{k+1}}\rightarrow \frac{F}{F_{k+1}}\rightarrow
\frac{F}{F_{k}}\rightarrow 1\text{,}
\end{equation*}%
and let $\left\{ E_{p,q}^{r}\right\} $ be the Hochschild-Serre spectral
sequence for the homology of this sequence. In particular, we have
\begin{equation*}
E_{p,q}^{2}=H_{p}\left( \frac{F}{F_{k}};H_{q}\left(
\frac{F_{k}}{F_{k+1}} \right) \right) .
\end{equation*}%
Then we have the differential
\begin{equation*}
d^{2}:E_{3,0}^{2}=H_{3}\left( \frac{F}{F_{k}}\right) \rightarrow
E_{1,1}^{2}=H_{1}\left( \frac{F}{F_{k}};H_{1}\left(
\frac{F_{k}}{F_{k+1}} \right) \right) \cong H_{1}\otimes
\frac{F_{k}}{F_{k+1}}\text{.}
\end{equation*}%
Finally, the refinement of Johnson's homomorphism is given by the
following theorem.

\begin{theorem}[Morita]
The composition $d^{2}\circ \tilde{\tau}_{k}$ coincides with
Johnson's homomorphism $\tau _{k}$ so that the following diagram
commutes.
\begin{equation*}
\begin{diagram}
               &                          & H_{3}\left( \frac{F}{F_{k}}\right) \\
               & \ruTo^{\tilde{\tau}_{k}} & \dTo_{d^{2}}                       \\
\mathcal{J}(k) & \rTo_{\tau _{k}}         & H_{1}\otimes \frac{F_{k}}{F_{k+1}} \\
\end{diagram}
\end{equation*}
\end{theorem}

\begin{theorem}[Morita]
Let $D_{k}(H_{1})$ be the submodule of $H_{1}\otimes
F_{k}/F_{k+1}$ defined to be the kernel of the natural surjection
\begin{equation*}
H_{1}\otimes \frac{F_{k}}{F_{k+1}}\rightarrow \frac{F_{k+1}}{F_{k+2}}
\end{equation*}%
given by the Lie bracket map $\left( w,\xi \right) \mapsto \left[
w,\xi \right] $. Then the image of the Johnson homomorphism $\tau
_{k}:\mathcal{J}(k)\rightarrow H_{1}\otimes F_{k}/F_{k+1}$ is
contained in $D_{k}(H_{1})$ so that we can write $\tau
_{k}:\mathcal{J}(k)\rightarrow D_{k}(H_{1}).$
\end{theorem}

A short remark about this theorem is perhaps in order. It is known
that the image of $\tau _{2}$ is exactly equal to $D_{2}(H_{1})$,
and the image of $\tau _{3}$ is a submodule of $D_{3}(H_{1})$ of
index a power of 2. Thus $\im\tau _{3}$ and $D_{3}(H_{1})$ have
the same rank. However, for $k\geq 4$, $k$ even, the rank of
$\im\tau _{k}$ is smaller than the rank of $D_{k}(H_{1})$. Please
see \cite{[Mo]} for more details.

\section{Birman-Craggs Homomorphism\label{(BCH)}}

As mentioned at the end of Section \ref{(JOD)} the Johnson
homomorphism $\tau _{2}$ only detects the free abelian part of the
abelianization of the Torelli group $\mathcal{J}(2),$ and some
2-torsion remains undetected. In this section we will say a word
about this 2-torsion. In \cite{[BC]} Birman and Craggs defined a
(finite) collection of abelian quotients of $\mathcal{J}(2)$ given
by homomorphisms onto $\mathbb{Z}_{2}.$ Here we will give a
description of these homomorphisms that is due to Johnson
\cite{[J1]}. This somewhat more tractable description is different
than (yet equivalent to) Birman and Craggs' original definition,
and it enabled Johnson to give the number of distinct
Birman-Craggs homomorphisms.

Consider the surface $\Sigma _{g,1},$ and let $f\in
\mathcal{J}(2)$. The definition of $\Gamma _{g,1}$ requires that
$f$\ be the identity on $\partial \Sigma _{g,1}.$ Thus $f$\ can
easily be extended to a homeomorphism of the closed surface
$\Sigma _{g}.$ Let $h:\Sigma _{g}\rightarrow S^{3}$ be a Heegaard
embedding of $\Sigma _{g}$ into the 3-sphere $S^{3}$, i.e. $\Sigma
_{g}$ bounds handlebodies on both sides in $S^{3}$. Now cut
$S^{3}$ open along $h(\Sigma _{g})$ and reglue the two pieces
using $f\in \mathcal{J}(2)$. The resulting manifold $S_{h,f}^{3}$
is a homology $S^{3}$, and its Rochlin invariant $\mu
(S_{h,f}^{3})\in \mathbb{Z}_{2}$ is defined.

In general, any closed, connected 3-manifold $M$, together with a
fixed trivialization of its tangent bundle over the 2-skeleton, is
the boundary of a 4-manifold $W$ whose tangent bundle can be
trivialized in a compatible fashion. If $s$ denotes the choice of
stable trivialization of the tangent bundle of $M$ over the
2-skeleton, then the \textit{Rochlin invariant} $\mu (M,s)\in
\mathbb{Z}_{16}$ is defined to be the signature $\sigma (W)$
reduced modulo 16. If $M$ happens to be a homology $S^{3}$ then
$s$ is unique and $\sigma (W)$ is divisible by 8. Thus $\mu
(S_{h,f}^{3})=\mu (S_{h,f}^{3},s)=$ $\sigma (W)$ can be considered
an element of $\mathbb{Z}_{2}.$ For a fixed Heegaard embedding
$h:\Sigma _{g}\rightarrow S^{3},$ the \textit{Birman-Craggs
homomorphism }$\rho _{h}:\mathcal{J}(2)\rightarrow \mathbb{Z}_{2}$
is defined by $\rho _{h}(f)=\mu (S_{h,f}^{3}).$

By relating the Birman-Craggs homomorphisms to a
$\mathbb{Z}_{2}$-quadratic form, Johnson was able to show the
dependence of $\rho _{h}$ on the embedding $h:\Sigma
_{g}\rightarrow S^{3}.$ The $\mathbb{Z}_{2}$-quadratic form
$q:H_{1}(\Sigma_{g};\mathbb{Z}_{2})\rightarrow \mathbb{Z}_{2}$ is
defined as follows. Let $\left\langle ~,~\right\rangle $ be the
Seifert linking form on $H_{1}(\Sigma _{g};\mathbb{Z}_{2})$
induced by $h:\Sigma _{g}\rightarrow S^{3}$ defined by letting
$\left\langle x,y\right\rangle $ be the linking number (modulo 2)
of $h(x)$ and $h(y)^{+}$ in $S^{3}$, where $h(y)^{+}$ is the
positive push-off of $h(y)$ in the normal direction determined by
the orientations of $h(\Sigma _{g})$ and $S^{3}.$ Define
$q(x)=\left\langle x,x\right\rangle $, then it is a
$\mathbb{Z}_{2}$-quadratic form on
$H_{1}(\Sigma_{g};\mathbb{Z}_{2})$\ induced by the embedding $h.$
Because it is a quadratic form, $q$ satisfies
$q(x+y)=q(x)+q(y)+x\cdot y,$ where $x\cdot y$ is the intersection
pairing of $H_{1}(\Sigma _{g};\mathbb{Z}_{2}).$ Let $\left\{
x_{i},y_{i}\right\} , 1\leq i\leq g,$ denote the standard basis
for $H_{1}(\Sigma _{g};\mathbb{Z}_{2})$, and the \textit{Arf
invariant }of $\Sigma _{g}$ with respect to $q$ is defined to be
\begin{equation*}
\arf(\Sigma _{g},q)=\sum_{i=1}^{g}q(x_{i})q(y_{i})\text{ (mod 2).}
\end{equation*}

Johnson's main results from \cite{[J1]} are as follows. Suppose
$h_{1},$ $h_{2}:\Sigma _{g}\rightarrow S^{3}$ are both Heegaard
embeddings of the surface $\Sigma _{g}.$

\begin{theorem}[Johnson]
The embeddings $h_{1}$ and $h_{2}$ induce the same mod 2
self-linking form if and only if the Birman-Craggs homomorphisms
$\rho _{h_{1}}$ and $\rho _{h_{2}}$ are equal.
\end{theorem}

Therefore the homomorphism $\rho _{h}:\mathcal{J}(2)\rightarrow
\mathbb{Z}_{2}$ only depends on the quadratic form $q$ induced by
$h$, and we replace the notation $\rho _{h}$ with $\rho _{q}$ to
emphasize this fact. Moreover, the $\mathbb{Z}_{2}$-quadratic
forms $q$ which are induced by a Heegaard embedding $h$ are
exactly those that satisfy $\arf(\Sigma _{g},q)=0.$ Thus we are
able to enumerate $\left\{ \rho _{q}\right\} .$

\begin{corollary}[Johnson]
There are precisely $2^{g-1}\left( 2^{g}+1\right) $ distinct
Birman-Craggs homomorphisms $\rho _{q}:\mathcal{J}(2)\rightarrow
\mathbb{Z}_{2}$.
\end{corollary}

Johnson also provided a means of computing $\rho _{q}$ in terms of
the Arf invariant.
\begin{enumerate}
\item[(1)] If $f\in \mathcal{J}(2)$ is a Dehn twist about a
bounding simple closed curve $C$, then
\begin{equation*}
\rho _{q}(f)=\arf(\Sigma ^{\prime },q|_{\Sigma ^{\prime }}),
\end{equation*}%
where $\Sigma ^{\prime }$ is a subsurface of $\Sigma _{g}$ bounded
by $C$.

\item[(2)] If $f\in \mathcal{J}(2)$ is a composition of Dehn twists about
cobounding curves $C_{1}$ and $C_{2},$ then
\begin{equation*}
\rho _{q}(f)=\left\{
\begin{array}{cc}
0 & \text{if }q(C_{1})=q(C_{2})=1 \\
\arf(\Sigma ^{\prime },q|_{\Sigma ^{\prime }}) & \text{if }
q(C_{1})=q(C_{2})=0
\end{array}%
\right.
\end{equation*}%
where $\Sigma ^{\prime }$ is a subsurface of $\Sigma _{g}$
cobounded by $C_{1}$ and $C_{2}.$
\end{enumerate}

For genus $g=2$ surfaces, the Torelli group $\mathcal{J}(2)$ is
generated by the collection of all Dehn twists about bounding
simple closed curves. For genus $g\geq 3$, $\mathcal{J}(2)$ is
generated by the collection of all Dehn twists about genus 1
cobounding pairs of simple closed curves, i.e. pairs of
non-bounding, disjoint, homologous simple closed curves that
together bound a genus 1 subsurface. Thus the list above is
sufficient for computing $\rho _{q}(f)$ for any $f\in
\mathcal{J}(2)$.

\section{Abelianization of the Torelli Group\label{(ATG)}}

We are now prepared to say something about the abelianization
\begin{equation*}
H_{1}(\mathcal{J}(2);\mathbb{Z)}\cong \frac{\mathcal{J}(2)}{\left[
\mathcal{J}(2),\mathcal{J}(2)\right] }
\end{equation*}%
of the Torelli group $\mathcal{J}(2).$ In fact, the main result of
Johnson in \cite{[J6]} is that the Johnson homomorphism $\tau
_{2}:\mathcal{J}(2)\rightarrow D_{2}(H_{1})$ and the totality of
the Birman-Craggs homomorphisms $\rho
_{q}:\mathcal{J}(2)\rightarrow \mathbb{Z}_{2}$ completely
determine $H_{1}(\mathcal{J}(2);\mathbb{Z)}$

On the one hand, we have the composition
\begin{equation*}
\begin{diagram}
\frac{\mathcal{J}(2)}{\left[
\mathcal{J}(2),\mathcal{J}(2)\right]} & \rTeXonto &
\frac{\mathcal{J}(2)}{\mathcal{J}(3)} & \rTeXto^{\cong} &
D_{2}(H_{1})
\end{diagram}
\end{equation*}%
where the first map is the projection given by the fact that
$\left[ \mathcal{J}(2),\mathcal{J}(2)\right] \subset
\mathcal{J}(3)$ and the second map is given by $\tau _{2}.$ After
we tensor with the rationals $\mathbb{Q}$, Johnson shows that we
obtain an isomorphism
\begin{equation*}
\begin{diagram}
\frac{\mathcal{J}(2)}{\left[
\mathcal{J}(2),\mathcal{J}(2)\right]} \otimes ~\mathbb{Q} &
\rTeXto^{\cong} & \frac{\mathcal{J}(2)}{\mathcal{J}(3)} \otimes
~\mathbb{Q} & \rTeXto^{\cong} & D_{2}(H_{1}) \otimes ~\mathbb{Q}
\end{diagram}.
\end{equation*}%
Thus we have $H_{1}(\mathcal{J}(2);\mathbb{Q)}\cong
\mathcal{J}(2)/\mathcal{J}(3)\otimes \mathbb{Q}.$

On the other hand, consider the totality of the Birman-Craggs
homomorphisms $\left\{ \rho _{q}\right\} $, and let
\begin{equation*}
\mathcal{C=}\bigcap\limits_{q}\ker \rho _{q}
\end{equation*}%
be the common kernel of all $\rho _{q}$ for all $q$ which satisfy
$\arf(\Sigma _{g},q)=0.$ Also let $\mathcal{J}(2)^{2}$ represent
the subgroup generated by all squares in $\mathcal{J}(2),$ and let
$\mathcal{O}_{q}$ be the subgroup of the mapping class group
$\Gamma _{g,1}$ which acts trivially on $H_{1}(\Sigma
_{g,1};\mathbb{Z}_{2}).$ That is, $\mathcal{O}_{q}$ consists of
those homeomorphisms which preserve the quadratic form $q.$ Then,
by using the theory of Boolean quadratic and cubic forms, Johnson
showed that
\begin{equation*}
\mathcal{C}=\mathcal{J}(2)^{2}=\left[ \mathcal{O}_{q},\mathcal{J}(2)\right] .
\end{equation*}%
Finally he showed that the commutator subgroup of $\mathcal{J}(2)$
is given by
\begin{equation*}
\left[ \mathcal{J}(2),\mathcal{J}(2)\right] =\mathcal{C}\cap \ker
\tau _{2}=\mathcal{C\cap J}(3).
\end{equation*}%
Thus we can completely determine
$H_{1}(\mathcal{J}(2);\mathbb{Z)}\cong \mathcal{J}(2)/\left[
\mathcal{J}(2),\mathcal{J}(2)\right] $ from the homomorphisms
$\left\{ \tau _{2},\rho _{q}\right\} .$

\section{A Bordism Representation of the Mapping Class
Group\label{(BRMCG)}}

\subsection{The Bordism Group $\Omega _{3}(X,A)$\label{(BG)}}

Let $\left( X,A\right) $ be a pair of topological spaces
$A\subseteq X$. The \textit{3-dimension oriented relative bordism
group} $\Omega _{3}(X,A)$ is defined to be the set of
\textit{bordism classes} of triples $\left( M,\partial M,\phi
\right) $ consisting of a compact, oriented 3-manifold $M$ with
boundary $\partial M$ and a continuous map $\phi :\left(
M,\partial M\right) \rightarrow \left( X,A\right) $. The triples
$\left( M_{0},\partial M_{0},\phi _{0}\right) $ and $\left(
M_{1},\partial M_{1},\phi _{1}\right) $ are equivalent, or
\textit{bordant over }$\left( X,A\right) $, if there exists a
triple $\left( W,\partial W,\Phi \right) $ consisting of a
compact, oriented 4-manifold $W$ with boundary $\partial W=\left(
M_{0}\amalg -M_{1}\right) \cup _{\partial M}M$ and a continuous
map $\Phi :\left( W,\partial W\right) \rightarrow \left(
X,A\right) $ satisfying $\Phi |_{M_{i}}=\phi _{i}$ and $\Phi
(M)\subset A.$ We also require that $\partial M=\partial
M_{0}\amalg -\partial M_{1}$ so that $\partial W$ is a closed
3-manifold.

\begin{figure}[h]
 \centering
 \includegraphics{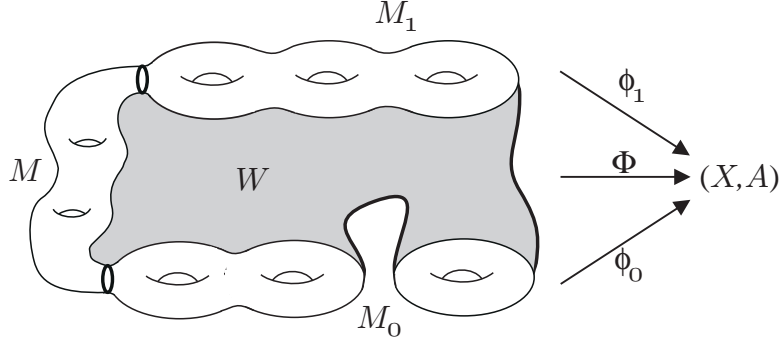}
 \caption{A relative bordism over $\left( X,A\right) $}
 \label{fig bordism def}
\end{figure}

A triple $\left( M,\partial M,\phi \right) $ is said to be
\textit{null-bordant} (or trivial) over $\left( X,A\right) $ if it
bounds $\left( W,\partial W,\Phi \right) ,$ that is, if it is
bordant to the empty set $\varnothing .$ The set $\Omega_{3}(X,A)$
forms a group with the operation of disjoint union and identity
element $\varnothing .$ In the case that $A=\varnothing $, we may
write $\Omega _{3}(X)=\Omega _{3}(X,\varnothing )$ and restrict
our definition to pairs $\left( M,\phi \right) =\left(
M,\varnothing ,\phi \right) $ of closed, oriented 3-manifolds.

\subsection{A Bordism Invariant of $\mathcal{J}(k)$\label{(BordInvt)}}

The purpose of this section is to analyze $\mathcal{J}(k)$ from
the point of view of bordism theory. Let $F=\pi_{1}(\Sigma_{g,1})$
as before, and consider the pair $\left( K(F/F_{k},1),\zeta
\right) $, where $K(F/F_{k},1)$ is an Eilenberg-MacLane space and
$\zeta \subset K(F/F_{k},1)$ is an $S^{1}$ corresponding to the
image of $\partial \Sigma _{g,1}$ under a continuous map $\Sigma
_{g,1}\rightarrow K(F/F_{k},1)$ induced by the canonical
projection $F\twoheadrightarrow F/F_{k}$. We denote the bordism
group over $\left( K(F/F_{k},1),\zeta \right) $ by $\Omega
_{3}(F/F_{k},\zeta ).$ Moreover, we have an isomorphism $j_{\ast
}:\Omega _{3}(F/F_{k})\rightarrow \Omega _{3}(F/F_{k},\zeta )$
induced by the inclusion map $j:\left( K(F/F_{k},1),\varnothing
\right) \rightarrow \left( K(F/F_{k},1),\zeta \right) .$ We will
make use of both of these groups in what follows, but our main
focus will be on the group $\Omega _{3}(F/F_{k}).$

Below, in Theorem \ref{(SigmaHomom)}, we define a homomorphism
$\sigma_{k}:\mathcal{J}(k)\rightarrow \Omega _{3}\left(F/F_{k}
\right)$ whose kernel is $\ker \sigma
_{k}=\mathcal{J}\left(2k-1\right).$ Thus the image of $\sigma
_{k}$ is $\mathcal{J}(k)/\mathcal{J}\left(2k-1\right).$ We already
saw that the image of the Johnson homomorphism $\tau _{k}$ is
$\mathcal{J}(k)/\mathcal{J}\left(k+1\right).$ However, since we
know $\left[\mathcal{J}(k),\mathcal{J}(k)\right] \subset
\mathcal{J}\left(2k-1\right)\subset \mathcal{J}\left(k+1\right),$
the image of this new homomorphism $\sigma _{k}$ is, in general,
much closer to the abelianization of $\mathcal{J}(k)$.

Consider a surface homeomorphism $f\in \mathcal{J}(k)$ for some
$k\geq 2$. As before let $T_{f,1}$ be the mapping torus of $f$,
i.e. $\Sigma _{g,1}\times \left[ 0,1\right] $ with $x\times
\left\{ 0\right\} $ glued to $f(x)\times \left\{ 1\right\} $. The
boundary\ $\partial T_{f,1}$ of $T_{f,1}$ is the torus $\partial
\Sigma _{g,1}\times S^{1},$ and the mapping torus $T_{f,1}$ is an
(oriented) homology $\Sigma _{g,1}\times S^{1}$. Fixing a basis
$\left\{ \alpha _{1},...,\alpha _{2g}\right\} $ for the free group
$F=\pi _{1}(\Sigma _{g,1})$ gives a presentation of $\pi
_{1}(T_{f,1})$:
\begin{equation*}
\pi _{1}(T_{f,1})=\left\langle \alpha _{1},...,\alpha _{2g},\gamma
\mid \left[ \alpha _{1},\gamma \right] f_{\ast }(\alpha
_{1})\alpha _{1}^{-1},..., \left[ \alpha _{2g},\gamma \right]
f_{\ast }(\alpha _{2g})\alpha _{2g}^{-1}\right\rangle
\end{equation*}%
where $\gamma $ represents a generator of $\pi _{1}(S^{1})$. We
now wish to obtain a closed 3-manifold from $T_{f,1}$ by filling
in its boundary. Let $T_{f}^{\gamma }=T_{f,1}^{\gamma }$ be the
result of performing a Dehn filling along a curve on $\partial
T_{f,1}$ represented by the homotopy class $\gamma .$ That is,
$T_{f}^{\gamma }$ is obtained by filling in the torus $\partial
T_{f,1}\simeq \partial \Sigma _{g,1}\times S^{1}$ with the solid
torus $\partial \Sigma _{g,1}\times D^{2}.$ Then we also have a
presentation for $\pi _{1}(T_{f}^{\gamma })$:
\begin{equation*}
\pi _{1}(T_{f}^{\gamma })=\left\langle \alpha _{1},...,\alpha
_{2g}\mid f_{\ast }(\alpha _{1})\alpha _{1}^{-1},...,f_{\ast
}(\alpha _{2g})\alpha _{2g}^{-1}\right\rangle
\end{equation*}%
Note that if $f$ is isotopic to the identity, then $T_{f}^{\gamma
}$ is homeomorphic to the connected sum of $2g$ $\left(
S^{1}\times S^{2}\right)$'s.

Now for all $m\leq k$ we can define $\phi _{f,m}:\left(
T_{f,1},\partial T_{f,1}\right) \rightarrow \left(
K(F/F_{m},1),\zeta \right) $ to be a continuous map induced by the
canonical epimorphism
\begin{equation*}
\pi _{1}(T_{f,1})\twoheadrightarrow
\frac{\pi_{1}(T_{f,1})}{\left\langle \gamma ,\left(
\pi_{1}(T_{f,1})\right) _{m}\right\rangle }\cong \frac{F}{F_{m}}
\end{equation*}%
where the isomorphism requires the fact that $f\in
\mathcal{J}(k)\subset \mathcal{J}(m)$ (see Lemma \ref{(phi)}
below.) Also, since we kill the homotopy class $\gamma $ in our
construction of $T_{f}^{\gamma }$, the map $\phi _{f,m}$ extends
to a continuous map $\phi _{f,m}^{\gamma }:T_{f}^{\gamma
}\rightarrow K(F/F_{m},1),$ and $\phi _{f,m}^{\gamma }$ induces
the canonical epimorphism
\begin{equation*}
\pi _{1}(T_{f}^{\gamma })\twoheadrightarrow \frac{\pi
_{1}(T_{f}^{\gamma })}{\left( \pi _{1}(T_{f}^{\gamma })\right)
_{m}}\cong \frac{F}{F_{m}}.
\end{equation*}%
Moreover, we have the following lemma.

\begin{lemma}
\label{(phi)}The following are equivalent:

\begin{enumerate}
\item[(a)] $f\in \mathcal{J}(m)$,

\item[(b)] $\frac{\pi _{1}(T_{f}^{\gamma })}{\left( \pi
_{1}(T_{f}^{\gamma })\right) _{m}}\cong \frac{F}{F_{m}}$ and
$\frac{\pi _{1}(T_{f,1})}{\left\langle \gamma ,\left( \pi
_{1}(T_{f,1})\right) _{m}\right\rangle } \cong \frac{F}{F_{m}}$,
and

\item[(c)] the continuous maps $\phi _{f,m}^{\gamma }$ and $\phi
_{f,m}$ exist as defined.
\end{enumerate}
\end{lemma}

\begin{proof}
This is an obvious fact, but we wish to emphasize it because of
the important role it will play later.

\begin{proofnosquare}[(a)$\iff $(b)]
If $f\in \mathcal{J}(m)$ then the relations $\left[ \alpha
_{i},\gamma \right] f_{\ast }(\alpha _{i})\alpha _{i}^{-1}$ in
$\pi _{1}(T_{f,1})$ become trivial modulo $\left\langle \gamma
,\left( \pi _{1}(T_{f,1})\right) _{m}\right\rangle $ since
$f_{\ast }$ acts as the identity on $F/F_{m},$ and we clearly have
a homomorphism (in fact, an isomorphism.) On the other hand, no
such homomorphism exists if $f\notin \mathcal{J}(m)$ because the
relations $\left[ \alpha _{i},\gamma \right] f_{\ast }(\alpha
_{i})\alpha _{i}^{-1}\equiv f_{\ast }(\alpha _{i})\alpha
_{i}^{-1}$ (mod $\gamma$) are certainly not trivial modulo $\left(
\pi _{1}(T_{f,1})\right) _{m}$.
\end{proofnosquare}

\begin{proofnosquare}[\textit{(b)}$\iff $\textit{(c)}]
It is a well-known property of Eilenberg-MacLane spaces that
continuous maps into them are in one-to-one correspondence with
homomorphisms into their fundamental group. (See \cite{[Wh]}
Theorem V.4.3.) Thus $\phi _{f,m}$ (and similarly for $\phi
_{f,m}^{\gamma }$) is defined if and only if the homomorphism
\begin{equation*}
\pi _{1}(T_{f,1})\twoheadrightarrow \frac{\pi
_{1}(T_{f,1})}{\left\langle \gamma ,\left( \pi
_{1}(T_{f,1})\right) _{m}\right\rangle }\cong \frac{F}{F_{m}}
\end{equation*}%
exists.
\end{proofnosquare}
\end{proof}

Let us now consider the pair $\left( T_{f}^{\gamma },\phi
_{f,k}^{\gamma }\right) \in \Omega _{3}(F/F_{k})$. We introduce a
new homomorphism giving a representation of $\mathcal{J}(k) $
which is very geometric in nature.

\begin{theorem}
\label{(SigmaHomom)}The map
\begin{equation*}
\sigma _{k}:\mathcal{J}(k)\rightarrow \Omega _{3}\left(
\frac{F}{F_{k}} \right)
\end{equation*}%
defined by $\sigma _{k}(f)=\left( T_{f}^{\gamma },\phi
_{f,k}^{\gamma }\right) $\ is a well-defined homomorphism.
\end{theorem}

We point out that one can similarly define a homomorphism into the
relative bordism group $\mathcal{J}(k)\rightarrow \Omega
_{3}(F/F_{k},\zeta ) $ which sends a mapping class $f\in
\mathcal{J}(k)$ to $\left( T_{f,1},\partial T_{f,1},\phi
_{f,k}\right) $. However, we will mainly focus on the homomorphism
given in Theorem \ref{(SigmaHomom)}.

\begin{proof}
Consider two homeomorphisms $f,g\in \mathcal{J}(k)$ for the
oriented surface $\Sigma _{1}$ with one boundary component. If $f$
and $g$ are isotopic, i.e. they represent the same mapping class,
then of course $T_{f}^{\gamma }$ and $T_{g}^{\gamma }$ are
homeomorphic and $\left( T_{f}^{\gamma },\phi _{f,k}^{\gamma
}\right) $ and $\left( T_{g}^{\gamma },\phi _{g,k}^{\gamma
}\right) $ are bordant. Thus $\sigma _{k}$ is certainly
well-defined.

To show $\sigma _{k}$ is indeed a homomorphism we need to show
that $\left( T_{f}^{\gamma },\phi _{f,k}^{\gamma }\right) \amalg
\left( T_{g}^{\gamma },\phi _{g,k}^{\gamma }\right) $ is bordant
to $\left( T_{f\circ g}^{\gamma },\phi _{f\circ g,k}^{\gamma
}\right) $ in $\Omega _{3}(F/F_{k})$ for any mapping classes
$f,g\in \mathcal{J}(k).$ To do so, we simply construct a bordism,
i.e. we build a 4-manifold $W$ and continuous map $\Phi
:W\rightarrow K(F/F_{k},1)$ with boundary given by
\begin{equation*}
\left( \partial W,\Phi |_{\partial W}\right) =\left[ \left(
T_{f}^{\gamma },\phi _{f,k}^{\gamma }\right) \amalg \left(
T_{g}^{\gamma },\phi _{g,k}^{\gamma }\right) \right] \amalg
-\left( T_{f\circ g}^{\gamma },\phi _{f\circ g,k}^{\gamma
}\right).
\end{equation*}

We begin by first constructing a 4-manifold between the mapping
tori $T_{f,1}\amalg T_{g,1}$ and $T_{f\circ g,1}.$ Recall that
\begin{equation*}
T_{f,1}=\frac{\Sigma _{1}\times \left[ 0,1\right] }{\left(
x,0\right) \sim \left( f(x),1\right) }.
\end{equation*}%
We may also consider $T_{f\circ g,1}$ in pieces as depicted in
Figure \ref{fig blocks}.
\begin{figure}[h]
 \centering
 \includegraphics{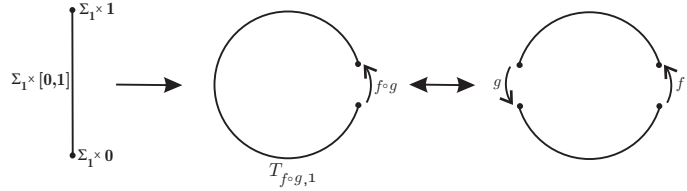}
 \caption{$T_{f\circ g,1}$ considered in pieces}
 \label{fig blocks}
\end{figure}

That is,
\begin{equation*}
T_{f\circ g,1}=\frac{\Sigma _{1}\times \left[ 0,1\right] }{\left(
x,0\right) \sim \left( f\left( g(x)\right) ,1\right) }\cong
\frac{\left( \Sigma _{1}\times \left[ 0,\frac{1}{2}\right] \right)
\cup \left( \Sigma _{1}\times \left[ \frac{1}{2},1\right] \right)
}{\left( x,0\right) \sim \left( f(x),1\right) ,\left(
x,\frac{1}{2}\right) \sim \left( g(x),\frac{1}{2} \right) }
\end{equation*}%

We can assume there is a product neighborhood of $\Sigma
_{1}\times \left\{ \frac{1}{2}\right\} $ in $T_{f,1}$, i.e. a
cylinder $\left( \Sigma _{1}\times \left\{ \frac{1}{2} \right\}
\right) \times \left[ -\varepsilon ,\varepsilon \right] .$ Let $
V=\left( T_{f,1}\amalg T_{g,1}\right) \times \left[ 0,1\right] .$
Then $V$ has boundary
\begin{equation*}
\partial V=\left( T_{f,1}\amalg T_{g,1}\right) \times \left\{ 0
\right\} \cup -\left( T_{f,1}\amalg T_{g,1}\right) \times \left\{ 1
\right\} \cup \left( \partial T_{f,1}\amalg \partial T_{g,1}
\right) \times \left[ 0,1\right] .
\end{equation*}%
Now consider the piece $\left( T_{f,1}\amalg T_{g,1}\right) \times
\left\{ 1\right\} $ of $\partial V$ and attach a 4-dimensional
\textquotedblleft strip\textquotedblright\ $\Sigma _{1}\times
\left[ -\varepsilon ,\varepsilon \right] \times \left[ -\delta
,\delta \right] $ to $\left( T_{f,1}\amalg T_{g,1}\right) \times
\left\{ 1\right\} $ by gluing $\Sigma _{1}\times \left[
-\varepsilon ,\varepsilon \right] \times \left\{ -\delta \right\}
$ to the neighborhood $\left( \Sigma _{1}\times \left\{
\frac{1}{2}\right\} \right) \times \left[ -\varepsilon
,\varepsilon \right] $ in $T_{f,1}$ and gluing $\Sigma _{1}\times
\left[ -\varepsilon ,\varepsilon \right] \times \left\{ \delta
\right\} $ to the neighborhood $\left( \Sigma _{1}\times \left\{
\frac{1}{2}\right\} \right) \times \left[ -\varepsilon
,\varepsilon \right] $ in $T_{g,1}.$ Let $V^{\prime }$ be the
result of this gluing, then
\begin{eqnarray*}
\partial V^{\prime } &=&\left( \left( T_{f,1}\amalg T_{g,1}
\right) \times \left\{ 0\right\} \right) \cup \left( -\left(
T_{f\circ g,1}\right) \times \left\{ 1\right\} \right) \cup
\left( \left( \partial T_{f,1}\amalg \partial T_{g,1}\right)
\times \left[ 0,1\right] \right) \\
&&\cup \left( \partial \Sigma _{1}\times \left[ -\varepsilon
,\varepsilon \right] \times \left[ -\delta ,\delta \right]
\right).
\end{eqnarray*}%
\begin{figure}[h]
 \centering
 \includegraphics{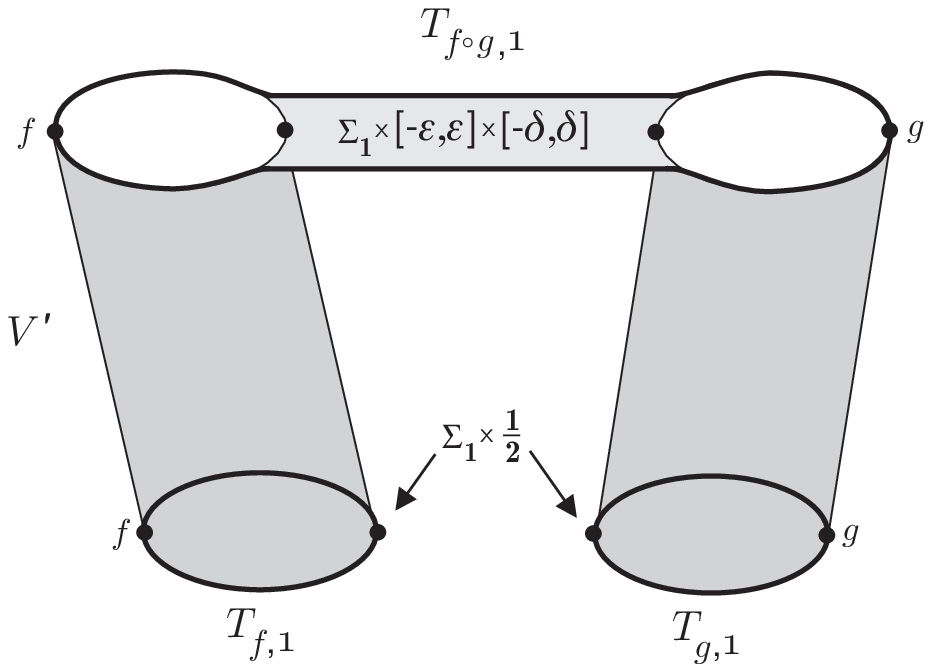}
 \caption{The 4-manifold $V^{\prime }$}
 \label{fig bordism}
\end{figure}

We now fill in the boundary component $\left( \partial
T_{f,1}\amalg \partial T_{g,1}\right) \times \left[ 0,1\right] $
with
\begin{equation}
\left( \left( \partial \Sigma _{1}\times D^{2}\right) \amalg
\left( \partial \Sigma _{1}\times D^{2}\right) \right) \times
\left[ 0,1\right] \tag{$\star $}  \label{(fill)}
\end{equation}%
to obtain a new 4-manifold $W$. At one end, this has the effect of
filling in the boundary of $\left( T_{f,1}\amalg T_{g,1}\right)
\times \left\{ 0\right\} ,$ thus creating $\left( T_{f}^{\gamma
}\amalg T_{g}^{\gamma }\right) \times \left\{ 0\right\} .$ At the
other end, we had already filled in \textit{some} of the boundary
of $T_{f\circ g,1}\times \left\{ 1\right\} $ with $\left(
\partial \Sigma _{1}\times \left[ -\varepsilon ,\varepsilon
\right] \times \left[ -\delta ,\delta \right] \right) $ above, and
the filling by (\ref{(fill)}) has the effect of filling in the
rest of the boundary of $T_{f\circ g,1}\times \left\{ 1\right\} .$
Thus we have actually created $T_{f\circ g}^{\gamma }\times
\left\{ 1\right\} .$ Therefore we have created a 4-manifold $W$
with boundary $\partial W=\left( T_{f}^{\gamma }\amalg
T_{g}^{\gamma }\right) \amalg -T_{f\circ g}^{\gamma }$. Also, the
continuous map $\phi _{f,k}\amalg \phi _{g,k}$ clearly extends
over $V=\left( T_{f,1}\amalg T_{g,1}\right) \times \left[
0,1\right] $. It is also easy to see that it extends over
$V^{\prime }$ as well since $\Sigma _{1}\times \left[ -\varepsilon
,\varepsilon \right] \times \left[ -\delta ,\delta \right] $
deformation retracts to $\Sigma _{1}.$ Finally it extends to a
continuous map $\Phi :W\rightarrow K(F/F_{k},1)$ in a similar way
that $\phi _{f,k}$ extends to $\phi _{f,k}^{\gamma }.$ Therefore
$\left( T_{f}^{\gamma },\phi _{f,k}^{\gamma }\right) \amalg \left(
T_{g}^{\gamma },\phi _{g,k}^{\gamma }\right) $ is bordant to
$\left( T_{f\circ g}^{\gamma },\phi _{f\circ g,k}^{\gamma }\right)
$ in $\Omega _{3}(F/F_{k}),$ and we have completed the proof of
Theorem \ref{(SigmaHomom)}.
\end{proof}

Notice that if a surface homeomorphism $f$ is isotopic to the
identity then its mapping class is in $\mathcal{J}(k)$ for all
$k,$ and $\left( T_{f,1},\partial T_{f,1},\phi _{f,k}\right)
=\left( T_{id,1},\partial T_{id,1},\phi _{id,k}\right) $ and
$\left( T_{f}^{\gamma },\phi _{f,k}^{\gamma }\right) =\left(
T_{id}^{\gamma },\phi _{id,k}^{\gamma }\right) $ are null-bordant
in $\Omega _{3}(F/F_{k},\zeta )$ and $\Omega _{3}(F/F_{k}),$
respectively, since they each bound $\Sigma _{g,1}\times D^{2}$
and the respective maps clearly extend. (The definition of
relative bordism requires an \textquotedblleft
extra\textquotedblright\ boundary piece so that the boundary of
the 4-manifold is a closed 3-manifold. In the case of $\left(
T_{id,1},\partial T_{id,1}\right) $, the extra piece is simply the
solid torus $\partial \Sigma _{g,1}\times D^{2}$ used to construct
$T_{id}^{\gamma }$.)

It seems logical to ask when $\left( T_{f,1},\partial T_{f,1},\phi
_{f,k}\right) \in \Omega _{3}(F/F_{k},\zeta )$ and $\left(
T_{f}^{\gamma },\phi _{f,k}^{\gamma }\right) \in \Omega
_{3}(F/F_{k})$ are null-bordant for more general $f\in
\mathcal{J}(k).$ That is, what is the kernel of $\sigma _{k}$?
This is answered by the following theorem.

\begin{theorem}
\label{(nullbordant)}$\left( T_{f,1},\partial T_{f,1},\phi
_{f,k}\right) \in \Omega _{3}(F/F_{k},\zeta )$ and $\left(
T_{f}^{\gamma },\phi _{f,k}^{\gamma }\right) \in \Omega
_{3}(F/F_{k})$ are trivial if and only if $f\in
\mathcal{J}\left(2k-1\right).$
\end{theorem}

\begin{corollary}
\label{(nullbordantCor)}The kernel of the homomorphism $\sigma
_{k}$ is $\mathcal{J}\left(2k-1\right).${\hfill {$\square $}}
\end{corollary}

We also have the following generalization of Theorem
\ref{(nullbordant)} which is a corollary to the proof of Theorem
\ref{(SigmaHomom)}.

\begin{corollary}
\label{(SigmaHomomCor)}Consider $f,g\in \mathcal{J}(k).$ Then the
following are equivalent:

\begin{enumerate}
\item[(a)] $f\circ g^{-1}\in \mathcal{J}\left(2k-1\right),$

\item[(b)] $\left( T_{f}^{\gamma },\phi _{f,k}^{\gamma }\right) $
is bordant to $\left( T_{g}^{\gamma },\phi _{g,k}^{\gamma }\right)
$ in $\Omega _{3}(F/F_{k}),$

\item[(c)] $\left( T_{f,1},\partial T_{f,1},\phi _{f,k}\right) $
is bordant to $\left( T_{g,1},\partial T_{g,1},\phi _{g,k}\right)
$ in $\Omega _{3}(F/F_{k},\zeta ).$
\end{enumerate}
\end{corollary}

\begin{proof}
Suppose we have $\left( T_{f}^{\gamma },\phi _{f,k}^{\gamma
}\right) =\left( T_{g}^{\gamma },\phi _{g,k}^{\gamma }\right) $ in
$\Omega _{3}(F/F_{k}).$ This is equivalent to having $\left(
T_{f}^{\gamma },\phi _{f,k}^{\gamma }\right) \amalg \left(
T_{g^{-1}}^{\gamma },\phi _{g^{-1},k}^{\gamma }\right) =\left(
T_{g}^{\gamma },\phi _{g,k}^{\gamma }\right) \amalg \left(
T_{g^{-1}}^{\gamma },\phi _{g^{-1},k}^{\gamma }\right) .$ However,
we showed in the proof of Theorem \ref{(SigmaHomom)} that this is
equivalent to $\left( T_{f\circ g^{-1}}^{\gamma },\phi _{f\circ
g^{-1},k}^{\gamma }\right) =\left( T_{g\circ g^{-1}}^{\gamma
},\phi _{g\circ g^{-1},k}^{\gamma }\right)$ in $\Omega
_{3}(F/F_{k}).$ The latter is just $\left( T_{id}^{\gamma },\phi
_{id,k}^{\gamma }\right) $, which is nullbordant. Thus Theorem
\ref{(nullbordant)} says that this is equivalent to $f\circ
g^{-1}\in \mathcal{J}\left(2k-1\right).$ The equivalence of (c) is
proved similarly.
\end{proof}

\begin{proof}[Proof of Theorem \protect\ref{(nullbordant)}]
We prove the theorem for the pair $\left( T_{f}^{\gamma },\phi
_{f,k}^{\gamma }\right) \in \Omega _{3}(F/F_{k}),$ and the proof
for the triple $\left( T_{f,1},\partial T_{f,1},\phi _{f,k}\right)
\in \Omega _{3}(F/F_{k},\zeta )$ is completely analogous. Suppose
$f\in \mathcal{J}(m)$, then for $l\leq m$ let $\pi
_{m,l}:K(F/F_{m},1)\rightarrow K(F/F_{l},1)$ be the projection map
such that $\phi _{f,l}^{\gamma }=\pi _{m,l}\circ \phi
_{f,m}^{\gamma }.$

\begin{proofnosquare}[$(\Longleftarrow )$]
Let us first suppose that $f\in \mathcal{J}\left(2k-1\right).$
Then the pair $\left( T_{f}^{\gamma },\phi _{f,2k-1}^{\gamma
}\right) $ is defined and is an element of $\Omega
_{3}(F/F_{2k-1}).$ The following lemma is due to K. Igusa and K.
Orr (\cite{[IO]}, Theorem 6.7.)
\end{proofnosquare}

\begin{lemma}[Igusa-Orr]
Let $\left( \pi _{m,k}\right) _{\ast }$ be the induced map on
$H_{3}$ and consider $x\in H_{3}(F/F_{m}).$ Then $x\in \ker \left(
\pi _{m,k}\right) _{\ast }$ if and only if $x\in \image\left( \pi
_{2k-1,m}\right) _{\ast }$ for $k\leq m\leq 2k-1.$ In particular,
the homomorphism $\left( \pi _{2k-1,k}\right) _{\ast
}:H_{3}\left(F/F_{2k-1}\right) \rightarrow
H_{3}\left(F/F_{k}\right)$ is trivial.
\end{lemma}

We have the following corollary.

\begin{corollary}
\label{(IO)}The homomorphism
\begin{equation*}
\left( \pi _{2k-1,k}\right) _{\ast }:\Omega _{3}\left(
\frac{F}{F_{2k-1}} \right) \rightarrow \Omega _{3}\left(
\frac{F}{F_{k}}\right)
\end{equation*}%
is trivial. Moreover, a bordism class is in $\ker \left( \pi
_{m,k}\right) _{\ast }$ if and only if it lies in the image of
$\left( \pi _{2k-1,m}\right) _{\ast }$ for $k\leq m\leq 2k-1.$
\end{corollary}

\begin{proof}
In general, $\Omega _{n}(X,A)$ is the $n$-dimensional bordism
group, and it is an extraordinary homology theory. Using the
Atiyah-Hirzebruch spectral sequence, (see G. Whitehead \cite{[Wh]}
for details,) one can express $\Omega _{n}(X,A)$ in terms of
ordinary homology with coefficient group $\Omega _{q}, $ where
$\Omega _{q}=\Omega _{q}(\cdot )$ is the bordism group of a single
point. In particular, $E_{p,q}^{2}\cong H_{p}(X,A;\Omega _{q})$
and the boundary operator is $d_{p,q}^{2}:E_{p,q}^{2}\rightarrow
E_{p-2,q+1}^{2}, $ and $\Omega _{n}(X,A)$ is built using
$H_{p}(X,A;\Omega _{q})$ with $p+q=n. $ Now $\Omega _{0}\cong
\mathbb{Z}$ and $\Omega _{1},$ $\Omega _{2},$ and $\Omega _{3}$
are all trivial. So in the case $n=3$ we have $\Omega
_{3}(X,A)\cong H_{3}(X,A;\Omega _{0})\cong H_{3}(X,A).$ In fact,
the isomorphism is given by $\left( M,\partial M,\phi \right)
\longmapsto \phi _{\ast }(\left[ M,\partial M\right] )$ where
$\left[ M,\partial M\right] $ denotes the fundamental class in
$H_{3}(M,\partial M).$ Of course it follows directly that $\Omega
_{3}(F/F_{k})\cong H_{3}(F/F_{k})$ (and $\Omega _{3}(F/F_{k},\zeta
)\cong H_{3}(F/F_{k},\zeta ),$) and we have the following
commutative diagram:
\begin{equation*}
\begin{diagram}
H_{3}\left( \frac{F}{F_{2k-1}}\right)       & \rTo^{\left( \pi_{2k-1,m}\right) _{\ast }} & H_{3}\left( \frac{F}{F_{m}}\right) & \rTo^{\left( \pi _{m,k}\right) _{\ast }} & H_{3}\left(\frac{F}{F_{k}}\right)             \\
\dTo^{\cong}                                &                                            & \dTo^{\cong} & & \dTo^{\cong}                                                                                                 \\
\Omega _{3}\left( \frac{F}{F_{2k-1}}\right) & \rTo^{\left( \pi_{2k-1,m}\right) _{\ast }} & \Omega _{3}\left( \frac{F}{F_{m}}\right) & \rTo^{\left( \pi _{m,k}\right) _{\ast }} & \Omega _{3}\left(\frac{F}{F_{k}}\right) \\
\end{diagram}
\end{equation*}%
Since the map $\left( \pi _{2k-1,k}\right) _{\ast }$ on $H_{3}$ is
the zero-homomorphism, the conclusion of the first part of the
corollary is proved. The proof of the latter part is also
immediate.
\end{proof}

The image of $\left( T_{f}^{\gamma },\phi _{f,2k-1}^{\gamma
}\right) $ under $\left( \pi _{2k-1,k}\right) _{\ast }:\Omega
_{3}(F/F_{2k-1})\rightarrow \Omega _{3}(F/F_{k})$ is
\begin{equation*}
\left( \pi _{2k-1,k}\right) _{\ast }(T_{f}^{\gamma },\phi
_{f,2k-1}^{\gamma })=\left( T_{f}^{\gamma },\pi _{2k-1,k}\circ
\phi _{f,2k-1}^{\gamma }\right) =\left( T_{f}^{\gamma },\phi
_{f,k}^{\gamma }\right) ,
\end{equation*}%
and Corollary \ref{(IO)} tells us that this image is trivial in
$\Omega _{3}(F/F_{k})$. Thus the condition $f\in
\mathcal{J}\left(2k-1\right)$ is certainly sufficient.

\begin{proofnosquare}[$(\implies )$]
The proof of the necessity of $f\in \mathcal{J}\left(2k-1\right)$
is much more subtle. If we assume that $\left( T_{f}^{\gamma
},\phi _{f,k}^{\gamma }\right) $ is trivial in $\Omega
_{3}(F/F_{k}),$ then Corollary \ref{(IO)} tells us that there is a
pair $\left( M,\phi \right) \in \Omega _{3}(F/F_{2k-1})$ that gets
sent to $\left( T_{f}^{\gamma },\phi _{f,k}^{\gamma }\right) $,
but we do not know anything more than that. We want to show that
$\phi _{f,2k-1}^{\gamma }$ is defined, and by Lemma \ref{(phi)} we
may achieve the desired conclusion $f\in
\mathcal{J}\left(2k-1\right)$.
\end{proofnosquare}

\begin{lemma}[Cochran-Gerges-Orr]
\label{(CGO6.8b)}Let $M$ be any oriented manifold such that $\pi
_{1}(M)=G$, and suppose $F$ is a free group. Then for $k>1$,
$G/G_{k}\cong F/F_{k}$ if and only if $H_{1}(M)$ is torsion-free
and all Massey products for $H^{1}(M)$ of length less than $k$
vanish. Under the latter conditions, any isomorphism
$G/G_{k-1}\cong F/F_{k-1}$ extends to $G/G_{k}\cong F/F_{k}.$
\end{lemma}

\begin{proof}
If $G/G_{k}\cong F/F_{k}$ for $k>1$, there is a continuous map
$\phi :M\rightarrow K(F/F_{k},1)$ that induces an isomorphism
$\phi ^{\ast }:H^{1}(F/F_{k})\rightarrow H^{1}(M)$ and $H_{1}(M)$
is clearly torsion-free. In \cite{[Or]} (Lemma 16) it is shown
that Massey products for $H^{1}(F/F_{k})$ of length less than $k$
vanish and length $k$ Massey products generate $H^{2}(F/F_{k})$.
Consider $x_{i}\in H^{1}(F/F_{k})$, then $\left\langle
x_{1},...,x_{n}\right\rangle =0$ for all $n<k.$ Also, the
naturality of Massey products (see property (\ref{(MP)}.2)) tells
us that $\phi ^{\ast }\left\langle x_{1},...,x_{n}\right\rangle
\subset \left\langle \phi ^{\ast }x_{1},...,\phi ^{\ast
}x_{n}\right\rangle .$ Thus for $n<k$ we certainly have $0\in
\left\langle \phi ^{\ast }x_{1},...,\phi ^{\ast
}x_{n}\right\rangle .$ However, the uniqueness of Massey products
given in property (\ref{(MP)}.1) tells us that the first nonzero
Massey product is uniquely defined, and we conclude that
$0=\left\langle \phi ^{\ast }x_{1},...,\phi ^{\ast
}x_{n}\right\rangle $ for $n<k.$ Therefore, since $\phi ^{\ast }$
is an isomorphism, all Massey products for $ H^{1}(M)$ of length
less than $k$ are zero.

On the other hand, if $H_{1}(M)$ is torsion-free and all Massey
products for $H^{1}(M)$ of length less than $k$ vanish then we
easily see that $H_{1}(M)\cong G/G_{2}\cong F/F_{2}.$ Now assume
by induction that $G/G_{k-1}\cong F/F_{k-1},$ and let $\psi
:F\rightarrow G$ be a homomorphism that induces this isomorphism.
We will extend this isomorphism to $G/G_{k}\cong F/F_{k}.$ It is
sufficient to show that $G_{k-1}/G_{k}\cong F_{k-1}/F_{k}.$ We
have the following commutative diagram
\begin{equation*}
\begin{diagram}
0                    & \rTo              & H_{2}\left( \frac{F}{F_{k-1}}\right) & \rTo^{\cong} & \frac{F_{k-1}}{F_{k}} & \rTo & 0 \\
                     &                   & \dTo^{\cong}_{\psi_{\ast}}           &              & \dTo_{\psi_{\ast}}    &      &   \\
H_{2}\left( G\right) & \rTo^{\pi_{\ast}} & H_{2}\left( \frac{G}{G_{k-1}}\right) & \rTo         & \frac{G_{k-1}}{G_{k}} & \rTo & 0 \\
\end{diagram}
\end{equation*}%
in which the horizontal maps are exact sequences. The fact that
the sequences are exact is a result of J. Stallings \cite{[St]}.
This diagram shows us that it is sufficient to show that $\pi
_{\ast }:H_{2}(G)\rightarrow H_{2}(G/G_{k-1})$ is trivial.
However, since $ H_{2}(M)\twoheadrightarrow H_{2}(G)$ is onto, we
need only show that $\pi _{\ast }:H_{2}(M)\rightarrow
H_{2}(G/G_{k-1})$ is trivial. As mentioned above, length $k-1$
Massey products $\left\langle x_{1},...,x_{k-1}\right\rangle $
generate $H^{2}(G/G_{k-1})\cong H^{2}(F/F_{k-1}).$ Then $\pi
^{\ast }\left\langle x_{1},...,x_{k-1}\right\rangle =\left\langle
\pi ^{\ast }x_{1},...,\pi ^{\ast }x_{k-1}\right\rangle =0$ since
length $k-1$ Massey products vanish for $M.$ Therefore $\pi ^{\ast
}$ and $\pi _{\ast }$ are trivial homomorphisms, and the
conclusion follows.
\end{proof}

A slightly more general version of the following lemma is proved
in \cite{[CGO]} (Theorem 4.2), and we include a proof here for
your convenience.

\begin{lemma}[Cochran-Gerges-Orr]
\label{(CGO4.2)}Suppose $M_{0}$ and  $M_{1}$ are closed, oriented
3-manifolds with $\pi _{1}(M_{0})=G_{0}$ and $\pi
_{1}(M_{1})=G_{1}$. Further suppose that there is an epimorphism
$\psi :G_{1}\rightarrow G_{0}/(G_{0})_{k}$, and then let $\phi
_{0}:M_{0}\rightarrow K(G_{0}/(G_{0})_{k},1)$ and $\phi
_{1}:M_{1}\rightarrow K(G_{0}/(G_{0})_{k},1)$ be continuous maps
so that $\left( \phi _{1}\right) _{\ast }=\psi $ and $\left(
M_{0},\phi _{0}\right) =\left( M_{1},\phi _{1}\right) $ in $\Omega
_{3}(G_{0}/(G_{0})_{k}).$ Then $\left( M_{0},\phi _{0}\right) $
and $\left( M_{1},\phi _{1}\right) $ are bordant over
$K(G_{0}/(G_{0})_{k},1)$ via a 4-manifold with only 2-handles (rel
$M_{0}$) whose attaching circles lie in $(G_{0})_{k}.$
\end{lemma}

\begin{proof}
Since $\left( M_{0},\phi _{0}\right) $ and $\left( M_{1},\phi
_{1}\right) $ are bordant, we know there exists a compact,
oriented 4-manifold $W$ and a continuous map $\Phi :W\rightarrow
K(G_{0}/(G_{0})_{k},1)$ such that $\partial (W,\Phi )=\left(
M_{0},\phi _{0}\right) \amalg \left( -M_{1},\phi _{1}\right) .$
$\Phi _{\ast }$ is already a surjection on $\pi _{1},$ and we can
make it an injection by performing surgery on loops in $W.$ Thus
we may assume $\Phi _{\ast }$ is an isomorphism. Now we choose a
handlebody structure for $W$ relative to $M_{0}$ with no 0-handles
or 4-handles. We then get rid of the 1-handles by trading them for
2-handles, i.e. we perform a surgery along a loop $c$ passing over
the 1-handles in the interior of $W$. In a similar manner, we can
get rid of the 3-handles by thinking of them as 1-handles relative
to $M_{1}.$ Let $V$ be the result of this handle swapping. We want
to make sure $\Phi $ extends to $V$, so because $\Phi _{\ast }$ is
an isomorphism it is necessary to make sure $c$ was null-homotopic
in $W$ since it is null-homotopic in $V.$ However, since $\left(
\phi _{0}\right) _{\ast }$ is surjective and $c$ is in the
interior of $W$, we can alter $c$ by a loop in $M_{0}$ so that the
altered $c$ is null-homotopic in $W$. Thus we may assume that the
2-handles are attached along loops $c$ in $(G_{0})_{k}.$
\end{proof}

\begin{lemma}
\label{(CGO6.11)}Let $M_{i}$ and $G_{i}$ $(i=0,1)$ be as in Lemma
\ref{(CGO4.2)}. For some free group $F$ suppose that $\phi
_{0}:M_{0}\rightarrow K(F/F_{k},1)$ and $\phi
_{1}:M_{1}\rightarrow K(F/F_{k},1)$ are continuous maps such that
$\phi _{0}$ induces an isomorphism $G_{0}/(G_{0})_{k}\cong
F/F_{k}$ and $\phi _{1}$ extends to a continuous map $\phi
_{1}:M_{1}\rightarrow K(F/F_{k+1},1)$ inducing
$G_{1}/(G_{1})_{k+1}\cong F/F_{k+1}.$ If $\left( M_{0},\phi
_{0}\right) $ is bordant to $\left( M_{1},\phi _{1}\right) $ in
$\Omega _{3}(F/F_{k}),$ then $\phi _{0}$ also extends so that it
induces $G_{0}/(G_{0})_{k+1}\cong F/F_{k+1}.$
\end{lemma}

\begin{proof}
Lemma \ref{(CGO4.2)} tells us there exists a bordism $\left(
W,\Phi \right) $ between $\left( M_{0},\phi _{0}\right) $ and
$\left( M_{1},\phi _{1}\right) $ over $K(F/F_{k},1)$ such that $W$
contains only 2-handles with attaching circles in $F_{k}$ and $\pi
_{1}(W)\cong F/F_{k}.$ Let $j_{i}:M_{i} \rightarrow W$ be
inclusion maps so that $\Phi \circ j_{i}=\phi _{i},$ $ i=0,1 $.
\begin{equation*}
\begin{diagram}
M_{i}        &                  &                          \\
\dTo^{j_{i}} & \rdTo^{\phi_{i}} &                          \\
 W           & \rTo_{\Phi}      & K\left( F/F_{k},1\right) \\
\end{diagram}
\end{equation*}%
Consider any collection $\left\{ x_{1},...,x_{k}\right\} \in
H^{1}(M_{0})$ of cohomology classes. Then choose $y_{i}\in
H^{1}(F/F_{k})$ so that $\phi _{0}^{\ast }(y_{i})=x_{i}.$ Since
$\pi _{1}(W)\cong F/F_{k}$ and $ G_{0}/(G_{0})_{k}\cong F/F_{k}$,
Lemma \ref{(CGO6.8b)} says that Massey products of length less
than $k$ vanish. Thus each of the following Massey products are
uniquely defined:
\begin{equation*}
\left\langle x_{1},...,x_{k}\right\rangle =\left\langle \phi
_{0}^{\ast }(y_{1}),...,\phi _{0}^{\ast }(y_{k})\right\rangle
=j_{0}^{\ast }\left\langle \Phi ^{\ast }(y_{1}),...,\Phi ^{\ast
}(y_{k})\right\rangle .
\end{equation*}%
If we can actually show that these Massey products vanish then we
can use Lemma \ref{(CGO6.8b)} to show that $\phi _{0}$ also
extends so as to induce $G_{0}/(G_{0})_{k+1}\cong F/F_{k+1},$ thus
completing the proof. We will show $\left\langle \Phi ^{\ast
}(y_{1}),...,\Phi ^{\ast }(y_{k})\right\rangle =0.$ Since
$G_{1}/(G_{1})_{k+1}\cong F/F_{k+1},$ Lemma \ref{(CGO6.8b)} says
Massey products for $H^{1}(M_{1})$ of length less than $k+1$
vanish. In particular, length $k$ Massey products are zero, thus
\begin{equation*}
j_{1}^{\ast }\left\langle \Phi ^{\ast }(y_{1}),...,\Phi ^{\ast
}(y_{k})\right\rangle =\left\langle \phi _{1}^{\ast
}(y_{1}),...,\phi _{1}^{\ast }(y_{k})\right\rangle =0.
\end{equation*}%
Now consider the following short exact sequence
\begin{equation*}
0\longrightarrow H_{2}(M_{1})\longrightarrow H_{2}(W)\longrightarrow
H_{2}(W,M_{1})\longrightarrow 0.
\end{equation*}%
Since we can view $W$ as $M_{1}\times \left[ 0,1\right] $ with
2-handles attached along circles in $F_{k},$ we see that
$H_{2}(W,M_{1})$ is a free abelian group generated by the cores of
the 2-handles (rel $M_{1}$). Thus this sequence splits and we can
write $H_{2}(W)\cong H_{2}(M_{1})\oplus H_{2}(W,M_{1}).$ Because
the attaching circles of the 2-handles lie in $F_{k} $, the images
of the generators of the latter summand are clearly spheres in
$K(F/F_{k},1).$ But since $K(F/F_{k},1)$ has trivial higher
homotopy groups, they must vanish in $H_{2}(F/F_{k}).$ Then by
considering the dual splitting $H^{2}(W)\cong H^{2}(M_{1})\oplus
H^{2}(W,M_{1})$ we know that the image of $H^{2}(F/F_{k})$ must be
contained in the summand $H^{2}(M_{1})$ of $H^{2}(W). $ Therefore
$j_{1}^{\ast }:H^{2}(W)\rightarrow H^{2}(M_{1})$ must be injective
on the image of $H^{2}(F/F_{k}),$ and we are able to conclude that
$\left\langle \Phi ^{\ast }(y_{1}),...,\Phi ^{\ast
}(y_{k})\right\rangle =0$.
\end{proof}

Consider the following result of V. Turaev \cite{[Tu]}.

\begin{lemma}[Turaev]
\label{(Turaev)}Let $G$ be a finitely generated nilpotent group of
nilpotency class at most $k\geq 1$, and let $\alpha \in H_{3}(G).$
Then there exists a closed, connected, oriented 3-manifold $M$ and
a continuous map $\psi :M\rightarrow K(G,1)$ such that $\psi
_{\ast }(\left[ M\right] )=\alpha $ and such that $\psi $ induces
an isomorphism $\pi _{1}(M)/(\pi _{1}(M))_{k}\cong G$ if and only
if
\begin{enumerate}
\item[(a)] the homomorphism $\tor (H^{2}(G))\rightarrow
\tor(H_{1}(G))$ defined by sending $x$ to $x\cap \alpha $ is an
isomorphism, and

\item[(b)] for any
$h\in H_{2}(G),$ there exists $y\in H^{1}(G)$ such that
\begin{equation*}
h-(y\cap \alpha )\in \ker \left( H_{2}(G)\rightarrow H_{2}\left(
\frac{G}{G_{k-1}}\right) \right) .
\end{equation*}
\end{enumerate}
\end{lemma}

\begin{corollary}
\label{(TuraevCor)}For any bordism class $\alpha \in \Omega
_{3}(F/F_{k})$ there exists a closed, connected, oriented
3-manifold $M$ and a continuous map $\psi :M\rightarrow
K(F/F_{k},1)$ such that $\left( M,\psi \right) =\alpha $ in
$\Omega _{3}(F/F_{k})$ and $\psi $ induces an isomorphism $\pi
_{1}(M)/(\pi _{1}(M))_{k}\cong F/F_{k}.$
\end{corollary}

\begin{proof}
We simply use the fact proved earlier that $\Omega
_{3}(F/F_{k})\cong H_{3}(F/F_{k})$ and apply the lemma in the case
that $G\cong F/F_{k}.$ The group $F/F_{k}$ is nilpotent with
nilpotency $k-1$. The groups $H^{2}(F/F_{k})$ and $H_{1}(F/F_{k})$
are each torsion-free. Thus condition (a) of Lemma \ref{(Turaev)}
is satisfied trivially. Using Stallings' exact sequence given in
\cite{[St]}, we have the following commutative diagram
\begin{equation*}
\begin{diagram}
H_{2}(F)=0 & \rTo & H_{2}\left( \frac{F}{F_{k}}\right)   & \rTo^{\cong} & \frac{F_{k}}{F_{k+1}} & \rTo & 0 \\
           &      & \dTo                                 &              & \dTo_{0-map}          &      &   \\
H_{2}(F)=0 & \rTo & H_{2}\left( \frac{F}{F_{k-1}}\right) & \rTo^{\cong} & \frac{F_{k-1}}{F_{k}} & \rTo & 0 \\
\end{diagram}
\end{equation*}%
which shows us that the map $H_{2}(F/F_{k})\rightarrow
H_{2}(F/F_{k-1})$ is the zero homomorphism. Thus condition (b) of
Lemma \ref{(Turaev)} is also satisfied trivially.
\end{proof}

\begin{lemma}
\label{(CGO6.12)}Let $M$ be any closed, oriented 3-manifold with
$\pi _{1}(M)=G$, and suppose there is a continuous map $\phi
_{k}:M\rightarrow K(F/F_{k},1)$ inducing an isomorphism
$G/G_{k}\cong F/F_{k}$ for some free group $F$. For $m\geq k,$
$\left( M,\phi _{k}\right) $ is in the image of $ \left( \pi
_{m,k}\right) _{\ast }:\Omega _{3}(F/F_{m})\rightarrow \Omega
_{3}(F/F_{k})$ if and only if the isomorphism $G/G_{k}\cong
F/F_{k}$ can be extended to an isomorphism $G/G_{m}\cong F/F_{m}$
induced by a continuous map $\phi _{m}:M\rightarrow K(F/F_{m},1)$
such that $\left( \pi _{m,k}\right) _{\ast }(M,\phi _{m})=\left(
M,\phi _{k}\right) .$
\end{lemma}

\begin{proof}
Suppose $\left( M,\phi _{k}\right) =\left( \pi _{m,k}\right)
_{\ast }(\alpha ),$ for some $\alpha \in \Omega _{3}(F/F_{m}).$ By
Corollary \ref{(TuraevCor)} there exists a closed, connected,
oriented 3-manifold $ M^{\prime }$ and a continuous map $\psi
:M^{\prime }\rightarrow K(F/F_{m},1)$ that induces an isomorphism
$\pi _{1}(M^{\prime })/(\pi _{1}(M^{\prime }))_{m}\cong F/F_{m}$
such that $\left( M^{\prime },\psi \right) =\alpha $ in $\Omega
_{3}(F/F_{m}).$ Therefore we have $\left( M,\phi _{k}\right)
=\left( \pi _{m,k}\right) _{\ast }(\alpha )=\left( \pi
_{m,k}\right) _{\ast }(M^{\prime },\psi )=\left( M^{\prime },\pi
_{m,k}\circ \psi \right) .$ Thus $ \left( M,\phi _{k}\right) $ and
$\left( M^{\prime },\pi _{m,k}\circ \psi \right) $ are bordant in
$\Omega _{3}(F/F_{k}).$ In the case $m=k+1$, Lemma \ref{(CGO6.11)}
gives the desired result. The case $m>k+1$ is achieved via
induction. The converse is clear.
\end{proof}

We are now ready to continue our proof of Theorem
\ref{(nullbordant)}. First, we are assuming that $\phi
_{f,k}^{\gamma }$ exists, so Lemma \ref{(phi)} tells us that at
the very least $f\in \mathcal{J}(k).$ We also assume that $\left(
T_{f}^{\gamma },\phi _{f,k}^{\gamma }\right) $ is trivial in
$\Omega _{3}(F/F_{k}).$ In particular, $\left( T_{f}^{\gamma
},\phi _{f,k}^{\gamma }\right) =\left( T_{id}^{\gamma },\phi
_{id,k}^{\gamma }\right) $ in $\Omega _{3}(F/F_{k}).$ Also, we
have
\begin{equation*}
\frac{\pi _{1}(T_{id}^{\gamma })}{\left( \pi _{1}(T_{id}^{\gamma
})\right) _{m}}\cong \frac{F}{F_{m}}\text{, for all }m\text{, and}
\end{equation*}%
\begin{equation*}
\frac{\pi _{1}(T_{f}^{\gamma })}{\left( \pi _{1}(T_{f}^{\gamma
})\right) _{m}}\cong \frac{F}{F_{m}}\text{, for all }m\leq k.
\end{equation*}%
Then by Lemma \ref{(CGO6.11)} we can extend the latter isomorphism
to
\begin{equation*}
\frac{\pi _{1}(T_{f}^{\gamma })}{\left( \pi _{1}(T_{f}^{\gamma
})\right) _{k+1}}\cong \frac{F}{F_{k+1}}.
\end{equation*}%
By Lemma \ref{(phi)} we are able to conclude that $f\in
\mathcal{J}\left(k+1\right)$ and that the continuous map $\phi
_{f,k+1}^{\gamma }$ exists, allowing us to consider $\left(
T_{f}^{\gamma },\phi _{f,k+1}^{\gamma }\right) \in \Omega
_{3}(F/F_{k+1}).$ Moreover, since we are assuming that $\left(
T_{f}^{\gamma },\phi _{f,k}^{\gamma }\right) $ is trivial in
$\Omega _{3}(F/F_{k}),$ we have
\begin{equation*}
\left( T_{f}^{\gamma },\phi _{f,k+1}^{\gamma }\right) \in \ker
\left( \Omega _{3}\left( \frac{F}{F_{k+1}}\right) \overset{\left(
\pi _{k+1,k}\right) _{\ast }}{\longrightarrow }\Omega _{3}\left(
\frac{F}{F_{k}}\right) \right) ,
\end{equation*}%
and by Corollary \ref{(IO)}
\begin{equation*}
\left( T_{f}^{\gamma },\phi _{f,k+1}^{\gamma }\right) \in
\image\left( \Omega _{3}\left( \frac{F}{F_{2k-1}}\right)
\overset{\left( \pi _{2k-1,k+1}\right) _{\ast }}{\longrightarrow
}\Omega _{3}\left( \frac{F}{F_{k+1}}\right) \right) .
\end{equation*}%
Thus Lemma \ref{(CGO6.12)} implies that the isomorphism
\begin{equation*}
\frac{\pi _{1}(T_{f}^{\gamma })}{\left( \pi _{1}(T_{f}^{\gamma
})\right) _{k+1}}\cong \frac{F}{F_{k+1}}
\end{equation*}%
extends to an isomorphism
\begin{equation*}
\frac{\pi _{1}(T_{f}^{\gamma })}{\left( \pi _{1}(T_{f}^{\gamma
})\right) _{2k-1}}\cong \frac{F}{F_{2k-1}}.
\end{equation*}%
Therefore, by Lemma \ref{(phi)}, we are able to conclude that
$f\in \mathcal{J}\left(2k-1\right).$ This completes the proof of
Theorem \ref{(nullbordant)}.
\end{proof}

\subsection{Relating $\protect\sigma _{k}$ to the Johnson
Homomorphism}

The goal of this section is to describe how the homomorphism
$\sigma _{k}: \mathcal{J}(k)\rightarrow \Omega _{3}(F/F_{k})$
relates to Johnson's homomorphism
$\tau_{k}:\mathcal{J}(k)\rightarrow D_{k}(H_{1})\subset
\Hom(H_{1},F_{k}/F_{k+1}) .$ It turns out that $ \tau _{k}$
factors through $\Omega _{3}(F/F_{k}).$ To see this, we will use
Kitano's definition of $\tau _{k}$ in terms of Massey products,
which we reviewed in Section \ref{(MPD)}.

Let $\mathfrak{X}$ denote the ring of formal power series in the
noncommutative variables $t_{1},...,t_{2g}$, and let
$\mathfrak{X}_{k}$ denote the submodule of $\mathfrak{X}$
corresponding to the degree $k$ part. Because $F_{k}/F_{k+1}$ is a
submodule of $\mathfrak{X}_{k}$, we can consider the homomorphism
\begin{equation*}
\tau _{k}:\mathcal{J}(k)\rightarrow \Hom(H_{1},\mathfrak{X}_{k})
\end{equation*}%
defined in Theorem \ref{(kitano)}. Recall from Section \ref{(MPD)}
that we are considering the following dual bases:
\begin{equation*}
\left\{ x_{1},...,x_{2g},y\right\} \in H_{1}(T_{f,1})\text{,}
\end{equation*}%
\begin{equation*}
\left\{ x_{1}^{\ast },...,x_{2g}^{\ast },y^{\ast }\right\} \in
H^{1}(T_{f,1}) \text{, and}
\end{equation*}%
\begin{equation*}
\left\{ X_{1},...,X_{2g}\right\} \in H_{2}(T_{f,1})\text{.}
\end{equation*}%
Define $\Psi ^{\prime }:\Omega _{3}(F/F_{k},\zeta )\rightarrow
\Hom(H_{1},\mathfrak{X}_{k})$ to be the map that sends the bordism
class $\left( T_{f,1},\partial T_{f,1},\phi _{f,k}\right) $ to the
homomorphism
\begin{equation*}
x_{i}\longmapsto \sum_{j_{1},...,j_{k}}\left\langle \left\langle
x_{j_{1}}^{\ast },...,x_{j_{k}}^{\ast }\right\rangle ,X_{i}\right\rangle
t_{j_{1}}\cdot \cdot \cdot t_{j_{k}}\text{.}
\end{equation*}%
Let $i_{\ast }:\Omega _{3}(F/F_{k})\rightarrow \Omega
_{3}(F/F_{k},\zeta )$ be the homomorphism induced by inclusion
which sends $\left( T_{f}^{\gamma },\phi _{f,k}^{\gamma }\right) $
to $\left( T_{f,1},\partial T_{f,1},\phi _{f,k}\right) $. Then we
define the homomorphism $\Psi :\Omega _{3}(F/F_{k})\rightarrow
\Hom(H_{1},\mathfrak{X}_{k})$ to be the composition $\Psi =\Psi
^{\prime }\circ i_{\ast }.$

\begin{theorem}
\label{(psi_well-def)}The map $\Psi $ is a well-defined
homomorphism. Moreover, the composition $\Psi \circ \sigma _{k}$
corresponds to the Johnson homomorphism $\tau _{k}$ so that we
have the following commutative
diagram.%
\begin{equation*}
\begin{diagram}
               &                    & \Omega _{3}\left( F/F_{k}\right) \\
               & \ruTo^{\sigma_{k}} & \dTo_{\Psi}                      \\
\mathcal{J}(k) & \rTo_{\tau_{k}}    & \Hom(H_{1},\mathfrak{X}_{k})     \\
\end{diagram}
\end{equation*}
\end{theorem}

\begin{proof}
We only need to show that $\Psi ^{\prime}:\Omega_{3}(F/F_{k},\zeta
)\rightarrow \Hom(H_{1},\mathfrak{X}_{k})$ is a well-defined
homomorphism, and the rest of the theorem clearly follows. We will
need the following lemma.

\begin{lemma}
\label{(CGO6.8a)}Suppose $\left( M_{0},\phi _{0}\right) $ and
$\left( M_{1},\phi _{1}\right) $ are closed, oriented 3-manifolds
with $\pi _{1}(M_{i})=G_{i}$ and continuous maps $\phi
_{i}:M_{i}\rightarrow K(G_{0}/(G_{0})_{k},1)$. Further suppose
$\phi _{1}$ induces an isomorphism $ G_{1}/(G_{1})_{k}\cong
G_{0}/(G_{0})_{k}.$ If $\left( M_{0},\phi _{0}\right) $ is bordant
to $\left( M_{1},\phi _{1}\right) $ in $\Omega
_{3}(G_{0}/(G_{0})_{k})$ and all Massey products for
$H^{1}(M_{0})$ of length less than $k$ vanish, then $\phi =\left(
\phi _{0}\right) _{\ast }^{-1}\circ \left( \phi _{1}\right) _{\ast
}:H_{1}(M_{1})\rightarrow H_{1}(M_{0})$ is an isomorphism such
that for $x_{i}\in H^{1}(M_{0}),$ $ \mathcal{E}_{i}\in
H_{2}(M_{0})$ Poincar\'{e} dual to $x_{i},$ and
$\mathcal{F}_{i}\in H_{2}(M_{1})$ Poincar\'{e} dual to $\phi
^{\ast }(x_{i})\in H^{1}(M_{1})$ we have
\begin{equation*}
\left\langle \left\langle x_{j_{1}},...,x_{j_{k}}\right\rangle
,\mathcal{E}_{i}\right\rangle =\left\langle \left\langle \phi
^{\ast }(x_{j_{1}}),...,\phi ^{\ast }(x_{j_{k}})\right\rangle
,\mathcal{F}_{i}\right\rangle
\end{equation*}%
where $\left\langle \text{ \ , \ }\right\rangle $ is the dual
pairing of $ H^{2}(M_{i})$ and $H_{2}(M_{i}).$
\end{lemma}

\begin{proof}
Since $\left( M_{0},\phi _{0}\right) $ is bordant to $\left(
M_{1},\phi _{1}\right) $ in $\Omega _{3}(G_{0}/(G_{0})_{k}),$ we
must also have $\left( \phi _{0}\right) _{\ast }(\left[
M_{0}\right] )=\left( \phi _{1}\right) _{\ast }(\left[
M_{1}\right] )$ in\ $H_{3}(G_{0}/(G_{0})_{k})$ where $\left[
M_{i}\right] $ is the fundamental class in $H_{3}(M_{i})$. The
bordism $\left( W,\Phi \right) $ between $\left( M_{0},\phi
_{0}\right) $ and $\left( M_{1},\phi _{1}\right) $ can be chosen
so that $\Phi $ induces an isomorphism $\pi _{1}(W)\cong
G_{0}/(G_{0})_{k}$ and the inclusion maps $j_{i}:M_{i}\rightarrow
W$ induce isomorphisms $G_{i}/\left( G_{i}\right) _{k}\cong \pi
_{1}(W)/\left( \pi _{1}(W)\right) _{k}.$

W. Dwyer proves in \cite{[Dw]} (Corollary 2.5) that for cohomology
classes $ \alpha _{i}\in H^{1}(W)$ we have $\left\langle \alpha
_{1},...,\alpha _{m}\right\rangle =0$ if and only if $j_{0}^{\ast
}\left\langle \alpha _{1},...,\alpha _{m}\right\rangle =0$ for
$m<k.$ However, by the naturality of Massey products given in
property (\ref{(MP)}.2), we know that $j_{0}^{\ast }\left\langle
\alpha _{1},...,\alpha _{m}\right\rangle \subset \left\langle
j_{0}^{\ast }(\alpha _{1}),...,j_{0}^{\ast }(\alpha
_{m})\right\rangle ,$ and the latter is $0$ since Massey products
of length less than $k$ vanish for $H^{1}(M_{0}).$ Thus
$\left\langle \alpha _{1},...,\alpha _{m}\right\rangle =0$ for all
$ \alpha _{i}\in H^{1}(W).$ Moreover, $j_{1}^{\ast
}:H^{1}(W)\rightarrow H^{1}(M_{1})$ is an isomorphism. Then for
any $y_{i}\in H^{1}(M_{1})$ there exists an $\alpha _{i}\in
H^{1}(W)$ such that $j_{1}^{\ast }(\alpha _{i})=y_{i}.$ Thus for
$m<k$ we have
\begin{eqnarray*}
\left\langle y_{1},...,y_{m}\right\rangle &=&\left\langle j_{1}^{\ast
}(\alpha _{1}),...,j_{1}^{\ast }(\alpha _{m})\right\rangle \\
&=&j_{1}^{\ast }\left\langle \alpha _{1},...,\alpha _{m}\right\rangle \\
&=&0
\end{eqnarray*}%
where the second equality follows from naturality. So then we have
shown that all Massey products of length less than $k$ vanish also
for $H^{1}(W)$ and $H^{1}(M_{1}).$ Thus Massey products for
$H^{1}(M_{0}),$ $H^{1}(M_{1}),$ and $H^{1}(W)$ of length $k$ are
uniquely defined.

Consider $x_{i}\in H^{1}(M_{0})$ with Poincar\'{e} dual
$\mathcal{E}_{i}\in H_{2}(M_{0}).$ Let $\mathcal{F}_{i}\in
H_{2}(M_{1})$ be Poincar\'{e} dual to $\phi ^{\ast }(x_{i})\in
H^{1}(M_{1}),$ where $\phi $ is the isomorphism given by the
composition $\phi =\left( \phi _{0}\right) _{\ast }^{-1}\circ
\left( \phi _{1}\right) _{\ast }:H_{1}(M_{1})\rightarrow
H_{1}(M_{0}).$ Then we have
\begin{eqnarray*}
\left( \phi _{0}\right) _{\ast }(\mathcal{E}_{i}) &=&\left( \phi
_{0}\right)_{\ast }(x_{i}\cap \left[ M_{0}\right] ) \\
&=&\left( \phi _{0}^{\ast }\right) ^{-1}(x_{i})\cap \left( \phi
_{0}\right)_{\ast }(\left[ M_{0}\right] ) \\
&=&\left( \left( \phi _{1}^{\ast }\right) ^{-1}\circ \phi ^{\ast
}\right) (x_{i})\cap \left( \phi _{1}\right) _{\ast }(\left[ M_{1}\right] ) \\
&=&\left( \phi _{1}\right) _{\ast }(\phi ^{\ast }(x_{i})\cap
\left[ M_{1}\right] ) \\
&=&\left( \phi _{1}\right) _{\ast }(\mathcal{F}_{i}),
\end{eqnarray*}%
where the second and fourth equalities follow from the naturality
of cap products.

Now choose $\beta _{i}\in H^{1}(G_{0}/(G_{0})_{k})$ such that
$\phi _{0}^{\ast }(\beta _{i})=x_{i}.$ Then
\begin{eqnarray*}
\left\langle \left\langle x_{j_{1}},...,x_{j_{k}}\right\rangle
,\mathcal{E} _{i}\right\rangle &=&\left\langle \left\langle \phi
_{0}^{\ast }(\beta _{j_{1}}),...,\phi _{0}^{\ast }(\beta
_{j_{k}})\right\rangle ,\mathcal{E}
_{i}\right\rangle \\
&=&\left\langle \left\langle \beta _{j_{1}},...,\beta
_{j_{k}}\right\rangle
,\left( \phi _{0}\right) _{\ast }(\mathcal{E}_{i})\right\rangle \\
&=&\left\langle \left\langle \beta _{j_{1}},...,\beta
_{j_{k}}\right\rangle
,\left( \phi _{1}\right) _{\ast }(\mathcal{F}_{i})\right\rangle \\
&=&\left\langle \left\langle \phi _{1}^{\ast }(\beta
_{j_{1}}),...,\phi _{1}^{\ast }(\beta _{j_{k}})\right\rangle ,\mathcal{F}_{i}\right\rangle \\
&=&\left\langle \left\langle \left( \phi ^{\ast }\circ \phi
_{0}^{\ast }\right) (\beta _{j_{1}}),...,\left( \phi ^{\ast }\circ
\phi _{0}^{\ast}\right) (\beta _{j_{k}})\right\rangle ,\mathcal{F}_{i}\right\rangle \\
&=&\left\langle \left\langle \phi ^{\ast }(x_{j_{1}}),...,\phi ^{\ast
}(x_{j_{k}})\right\rangle ,\mathcal{F}_{i}\right\rangle .
\end{eqnarray*}%
This completes the proof of the lemma.
\end{proof}

Consider the mapping classes $f,h\in \mathcal{J}(k).$ We have the
dual bases mentioned above for specific homology and cohomology
groups of $T_{f,1}.$ Consider the following dual bases defined in
the same manner for $T_{h,1}$:
\begin{equation*}
\left\{ w_{1},...,w_{2g},z\right\} \in H_{1}(T_{h,1})\text{,}
\end{equation*}%
\begin{equation*}
\left\{ w_{1}^{\ast },...,w_{2g}^{\ast },z^{\ast }\right\} \in
H^{1}(T_{h,1}) \text{, and}
\end{equation*}%
\begin{equation*}
\left\{ W_{1},...,W_{2g}\right\} \in H_{2}(T_{h,1}).
\end{equation*}

Recall that $T_{f}^{\gamma }$ was constructed from $T_{f,1}$ by
filling the boundary $\partial T_{f,1}=\partial \Sigma
_{g,1}\times S^{1}$ with a solid torus $\partial \Sigma
_{g,1}\times D^{2}.$ Let $\psi _{f}:T_{f,1}\rightarrow
T_{f}^{\gamma }$ be the inclusion map, and then we have a basis
\begin{equation*}
\left\{ a_{1}^{\ast },...,a_{2g}^{\ast }\right\} \in H^{1}(T_{f}^{\gamma })
\end{equation*}%
where $a_{i}^{\ast }=\left( \left( \psi _{f}\right) _{\ast
}(x_{i})\right) ^{\ast }.$ Since $x_{i}^{\ast }$ is the Hom dual
of $x_{i},$ by definition we have that the dual pairing is
$\left\langle x_{i}^{\ast },x_{j}\right\rangle =\delta _{ij}.$
Similarly $a_{i}^{\ast }$ is the Hom dual of $\left( \psi
_{f}\right) _{\ast }(x_{i}),$ and thus $\left\langle \psi
_{f}^{\ast }(a_{i}^{\ast }),x_{j}\right\rangle =\left\langle
a_{i}^{\ast },\left( \psi _{f}\right) _{\ast }(x_{i})\right\rangle
=\delta _{ij}.$ Note that this implies that $\psi _{f}^{\ast
}(a_{i}^{\ast })=x_{i}^{\ast }.$ Letting $A_{i}$ denote the
Poincar\'{e} dual of $ a_{i}^{\ast }$ gives a basis for
$H_{2}(T_{f}^{\gamma })$:
\begin{equation*}
\left\{ A_{1},...,A_{2g}\right\} \in H_{2}(T_{f}^{\gamma }).
\end{equation*}%
By carefully examining the following commutative diagram, we see
$\left( \psi _{f}\right) _{\ast }(X_{i})=A_{i}.$
\begin{equation*}
\begin{diagram}
H_{1}\left (T_{f,1},\partial T_{f,1}\right) & \rTo^{\cong}_{\Hom~dual} & H^{1}\left (T_{f,1},\partial T_{f,1}\right) & \rTo^{\cap \left[ T_{f,1},\partial T_{f,1}\right]} & H_{2}\left (T_{f,1} \right)        \\
\uTo^{j_{*}}                                &                          & \dTo^{j^{*}}                                &                                                    &                                    \\
H_{1}\left (T_{f,1}\right)                  & \rTo^{\cong}_{\Hom~dual} & H^{1}\left(T_{f,1}\right)                   &                                                    &\dTo^{\left(\psi_{f}\right)_{\ast}} \\
\dTo^{\left(\psi_{f}\right)_{\ast}}         &                          & \uTo^{\psi_{f}^{\ast}}                      &                                                    &                                    \\
H_{1}\left( T_{f}^{\gamma }\right)          & \rTo^{\cong}_{\Hom~dual} & H^{1}\left( T_{f}^{\gamma }\right)          & \rTo^{\cap \left[T_{f}^{\gamma}\right]}            & H_{2}\left (T_{f}^{\gamma}\right)  \\
\end{diagram}
\end{equation*}%
Finally, let $\bar{A}_{i}\in $ $H_{2}(T_{h}^{\gamma })$ denote the
Poincar\'{e} dual to $\phi ^{\ast }(a_{i}^{\ast }),$ where $\phi $
is the isomorphism guaranteed by the following corollary. Then for
$T_{h,1}$ we similarly have $\psi _{h}^{\ast }(\phi ^{\ast
}(a_{i}^{\ast }))=w_{i}^{\ast } $ and $\left( \psi _{h}\right)
_{\ast }(W_{i})=\bar{A}_{i}.$ We have the following immediate
corollary to Lemma \ref{(CGO6.8a)}.

\begin{corollary}
\label{(CGO6.8aCor)}If $\left( T_{f}^{\gamma },\phi _{f,k}^{\gamma
}\right) =\left( T_{h}^{\gamma },\phi _{h,k}^{\gamma }\right) $ in
$\Omega _{3}(F/F_{k}),$ then the isomorphism $\phi =( \phi _{0})
_{\ast }^{-1}\circ ( \phi _{1}) _{\ast }:H_{1}(T_{g}^{\gamma
})\rightarrow H_{1}(T_{f}^{\gamma })$ satisfies
\begin{equation*}
\left\langle \left\langle a_{j_{1}}^{\ast },...,a_{j_{k}}^{\ast
}\right\rangle ,A_{i}\right\rangle =\left\langle \left\langle \phi
^{\ast }(a_{j_{1}}^{\ast }),...,\phi ^{\ast }(a_{j_{k}}^{\ast
})\right\rangle ,\bar{A}_{i}\right\rangle
\end{equation*}
where $\bar{A}_{i}$ is Poincar\'{e} dual to $\phi ^{\ast
}(a_{i}^{\ast }).${\hfill {$\square $}}
\end{corollary}

\begin{lemma}
If $\left( T_{f,1},\partial T_{f,1},\phi _{f,k}\right) =\left(
T_{h,1},\partial T_{h,1},\phi _{h,k}\right) $ in $\Omega
_{3}(F/F_{k},\zeta ),$ then
\begin{equation*}
\left\langle \left\langle x_{j_{1}}^{\ast },...,x_{j_{k}}^{\ast
}\right\rangle ,X_{i}\right\rangle =\left\langle \left\langle
w_{j_{1}}^{\ast },...,w_{j_{k}}^{\ast }\right\rangle
,W_{i}\right\rangle .
\end{equation*}
\end{lemma}

\begin{proof}
Since $f,h\in \mathcal{J}(k)$, Theorem \ref{(kitcor)} says that
the Massey products of length less than $k$ for $\left(
T_{f,1},\partial T_{f,1}\right) $ and $\left( T_{h,1},\partial
T_{h,1}\right) $ must vanish. Thus $ \left\langle x_{j_{1}}^{\ast
},...,x_{j_{k}}^{\ast }\right\rangle $ and $ \left\langle
w_{j_{1}}^{\ast },...,w_{j_{k}}^{\ast }\right\rangle $ are
uniquely defined. By Corollary \ref{(SigmaHomomCor)}, we know that
$\left( T_{f}^{\gamma },\phi _{f,k}^{\gamma }\right) =\left(
T_{h}^{\gamma },\phi _{h,k}^{\gamma }\right) $ in $\Omega
_{3}(F/F_{k}).$ So we let $\phi $ be the isomorphism guaranteed by
Corollary \ref{(CGO6.8aCor)}. Then we have
\begin{eqnarray*}
\left\langle \left\langle x_{j_{1}}^{\ast },...,x_{j_{k}}^{\ast
}\right\rangle ,X_{i}\right\rangle &=&\left\langle \left\langle \psi
_{f}^{\ast }(a_{j_{1}}^{\ast }),...,\psi _{f}^{\ast }(a_{j_{k}}^{\ast
})\right\rangle ,X_{i}\right\rangle \\
&=&\left\langle \left\langle a_{j_{1}}^{\ast },...,a_{j_{k}}^{\ast
}\right\rangle ,\left( \psi _{f}\right) _{\ast }(X_{i})\right\rangle \\
&=&\left\langle \left\langle a_{j_{1}}^{\ast },...,a_{j_{k}}^{\ast
}\right\rangle ,A_{i}\right\rangle \\
&=&\left\langle \left\langle \phi ^{\ast }(a_{j_{1}}^{\ast }),...,\phi
^{\ast }(a_{j_{k}}^{\ast })\right\rangle ,\bar{A}_{i}\right\rangle \\
&=&\left\langle \left\langle \phi ^{\ast }(a_{j_{1}}^{\ast }),...,\phi
^{\ast }(a_{j_{k}}^{\ast })\right\rangle ,\left( \psi _{h}\right) _{\ast
}(W_{i})\right\rangle \\
&=&\left\langle \psi _{h}^{\ast }\left\langle \phi ^{\ast }(a_{j_{1}}^{\ast
}),...,\phi ^{\ast }(a_{j_{k}}^{\ast })\right\rangle ,W_{i}\right\rangle \\
&=&\left\langle \left\langle \psi _{h}^{\ast }(\phi ^{\ast }(a_{j_{1}}^{\ast
})),...,\psi _{h}^{\ast }(\phi ^{\ast }(a_{j_{k}}^{\ast }))\right\rangle
,W_{i}\right\rangle \\
&=&\left\langle \left\langle w_{j_{1}}^{\ast },...,w_{j_{k}}^{\ast
}\right\rangle ,W_{i}\right\rangle .
\end{eqnarray*}
\end{proof}

This proves that $\Psi ^{\prime }:\Omega _{3}(F/F_{k},\zeta
)\rightarrow \Hom(H_{1},\mathfrak{X}_{k})$ is a well-defined
homomorphism and completes the proof of Theorem
\ref{(psi_well-def)}.
\end{proof}

\subsection{Relating $\protect\sigma _{k}$ to Morita's Homomorphism}

As we have already seen in the proof of Corollary \ref{(IO)},
there is an isomorphism $\Phi :\Omega _{3}(F/F_{k})\rightarrow
H_{3}(F/F_{k})$ given by $\left( M,\phi \right) \mapsto \phi
_{\ast }(\left[ M\right] )$, where $ \left[ M\right] $ is the
fundamental class in $H_{3}(M).$ Because of this, one may guess
that there is a relationship between $\sigma _{k}:\mathcal{J}
(k)\rightarrow \Omega _{3}(F/F_{k})$ and Morita's refinement
$\tilde{\tau}_{k}:\mathcal{J}(k)\rightarrow H_{3}(F/F_{k})$
discussed in Section \ref{(MR)}. This assumption turns out to be
correct, and the two homomorphisms are in fact equivalent.
However, $\sigma _{k}$ gives a representation that is much more
geometric, and as we will see in Section \ref{(SBRMCG)}, $\sigma
_{k}$ leads to interesting questions that $\tilde{\tau}_{k}$ does
not.

\begin{theorem}
\label{(sigma_is_morita)}The homomorphism $\sigma _{k}:\mathcal{J}
(k)\rightarrow \Omega _{3}(F/F_{k})$ coincides with the Morita
refinement of the Johnson homomorphism so that we have a
commutative diagram.
\begin{equation*}
\begin{diagram}
               &                         & \Omega _{3}\left( \frac{F}{F_{k}}\right) \\
               & \ruTo^{\sigma_{k}}      & \dTo^{\cong}_{\Phi}                      \\
\mathcal{J}(k) & \rTo_{\tilde{\tau}_{k}} & H_{3}\left( \frac{F}{F_{k}}\right)       \\
\end{diagram}
\end{equation*}
\end{theorem}

\begin{corollary}
\label{(morita kernel)}The kernel of Morita's refinement
$\tilde{\tau}_{k}$ is $\mathcal{J}\left(2k-1\right).$
\end{corollary}

\begin{proof}
This is an immediate consequence of Theorem
\ref{(sigma_is_morita)} and Corollary \ref{(nullbordantCor)}.
\end{proof}

\begin{proof}[Proof of Theorem \protect\ref{(sigma_is_morita)}]
Consider a genus $g$ surface with one boundary component $\Sigma
=\Sigma _{g,1}$ and $f\in \mathcal{J}(k)$. Let $r:\Sigma \times
\left[ 0,1\right] \rightarrow \Sigma $ be a retraction, $\psi
:\Sigma \rightarrow K(F/F_{k},1)$ be a continuous map that induces
the canonical epimorphism $F\twoheadrightarrow F/F_{k},$ and $
i:K(F/F_{k},1)\rightarrow \left( K(F/F_{k},1),\zeta \right) $ be
the inclusion map. Also let $G:\Sigma \times \left[ 0,1\right]
\rightarrow \left( T_{f,1},\partial T_{f,1}\right) $ be the
composition of the \textquotedblleft gluing map\textquotedblright\
$\Sigma \times \left[ 0,1 \right] \rightarrow T_{f,1}$ and the
inclusion $T_{f,1}\rightarrow \left( T_{f,1},\partial
T_{f,1}\right) $. Recall that the maps $\phi _{f,k}$ and $ \phi
_{f,k}^{\gamma }$ defined at the beginning of Section
\ref{(BordInvt)} are defined only up to homotopy. We choose them
so that the following diagram commutes.
\begin{equation*}
\begin{diagram}
\Sigma\times\left[ 0,1\right] & \rTo^r & \Sigma            &                            &                                              \\
                              &        &                   &\rdTo^{\psi}                &                                              \\
\dTo^G                        &        & T_{f}^{\gamma}    & \rTo^{\phi_{f,k}^{\gamma}} & K\left( F/F_{k},1\right)                     \\
                              &        &                   &                            & \dTo_i                                       \\
(T_{f,1},\partial T_{f,1})    &        & \rTo^{\phi_{f,k}} &                            & \left(K\left (F/F_{k},1\right),\zeta \right) \\
\end{diagram}
\end{equation*}%
That is, we have $\phi _{f,k}\circ G=i\circ \psi \circ r.$

Consider the fundamental class $\left[ T_{f,1},\partial
T_{f,1}\right] \in H_{3}(T_{f,1},\partial T_{f,1}),$ and suppose
that $\left( t_{f},\partial t_{f}\right) \in
C_{3}(T_{f,1},\partial T_{f,1})$ is a corresponding relative
3-cycle. Now we choose a 2-chain $\sigma \in C_{2}(\Sigma \times
\left[ 0,1\right] )$ so that $\partial \sigma $ is in the homotopy
class of a simple closed curve on $\Sigma \times \left\{ 0\right\}
$ parallel to the boundary $\partial \Sigma \times \left\{
0\right\} .$ Let $\sigma $ also denote $r_{\#}(\sigma )\in
C_{2}(\Sigma ),$ and choose a 3-chain $\rho \in C_{3}(\Sigma
\times \left[ 0,1\right] )$ so that $G_{\#}(\rho )=\left(
t_{f},\partial t_{f}\right) $ and $\partial \rho =\sigma
-f_{\#}(\sigma )+\left( \partial \sigma \times \left[ 0,1\right]
\right) .$

Consider the restriction $r|_{\partial \Sigma \times \left[
0,1\right] }.$ Then $r_{\#}(\partial \sigma \times \left[
0,1\right] )=\varepsilon \in C_{2}(\partial \Sigma ),$ and
\begin{eqnarray*}
\partial r_{\#}(\rho ) &=&r_{\#}\partial (\rho ) \\
&=&r_{\#}(\sigma -f_{\#}(\sigma )+\left( \partial \sigma \times
\left[ 0,1 \right] \right) ) \\
&=&r_{\#}(\sigma )-r_{\#}(f_{\#}(\sigma ))+r_{\#}(\partial \sigma \times
\left[ 0,1\right] ) \\
&=&\sigma -f_{\#}(\sigma )+\varepsilon
\end{eqnarray*}%
Since $f$ is the identity on the boundary, we must have $\partial
\sigma -f_{\#}(\partial \sigma )=0,$ and therefore $0=\partial
(\partial r_{\#}(\rho ))=\partial (\sigma -f_{\#}(\sigma
)+\varepsilon )=\partial \sigma -f_{\#}(\partial \sigma )+\partial
\varepsilon =\partial \varepsilon. $ Since $H_{2}(\partial \Sigma
)$ is trivial, there must be a 3-chain $ \eta \in C_{3}(\partial
\Sigma )$ such that $\partial \eta =\varepsilon .$ Let $j:\Sigma
\rightarrow \Sigma \times \left[ 0,1\right] $ be the inclusion
map, and consider $j_{\#}(\eta )\in C_{3}(\partial \Sigma \times
\left[ 0,1 \right] )\rightarrow C_{3}(\Sigma \times \left[
0,1\right] ).$ Define $ c_{f}\in C_{3}(\Sigma )$ to be
\begin{eqnarray*}
c_{f} &=&r_{\#}(\rho -j_{\#}(\eta )) \\
&=&r_{\#}(\rho )-r_{\#}j_{\#}(\eta ) \\
&=&r_{\#}(\rho )-\eta .
\end{eqnarray*}%
Then $\partial c_{f}=\partial r_{\#}(\rho )-\partial \eta =\left(
\sigma -f_{\#}(\sigma )+\varepsilon \right) -\varepsilon =\sigma
-f_{\#}(\sigma ).$

Also, since $j_{\#}(\eta )\in C_{3}(\Sigma \times \left[
0,1\right] )$ is carried by the subcomplex $\partial \Sigma \times
\left[ 0,1\right] ,$ $ G_{\#}(j_{\#}(\eta ))$ must be carried by
$\partial T_{f,1}.$ Thus $ G_{\#}(j_{\#}(\eta ))=0,$ and
\begin{eqnarray*}
G_{\#}(\rho -j_{\#}(\eta )) &=&G_{\#}(\rho )-G_{\#}(j_{\#}(\eta )) \\
&=&\left( t_{f},\partial t_{f}\right) .
\end{eqnarray*}%
Let $\bar{c}_{f}=\psi _{\#}(c_{f})\in C_{3}(F/F_{k}).$ Then
$\bar{c}_{f}$ is a 3-cycle since $f\in \mathcal{J}(k)$ induces the
identity on $F/F_{k}$. Let $\left[ \bar{c}_{f}\right] \in
H_{3}(F/F_{k})$ denote the corresponding homology class, and
\begin{eqnarray*}
i_{\ast }(\left[ \bar{c}_{f}\right] ) &=&\left[ i_{\#}(\bar{c}_{f})\right] \\
&=&\left[ \left( i\circ \psi \right) _{\#}(c_{f})\right] \\
&=&\left[ \left( i\circ \psi \circ r\right) _{\#}(\rho -j_{\#}(\eta ))\right]
\\
&=&\left[ \left( \phi _{f,k}\circ G\right) _{\#}(\rho -j_{\#}(\eta ))\right]
\\
&=&\left[ \left( \phi _{f,k}\right) _{\#}(t_{f},\partial t_{f})\right] \\
&=&\left( \phi _{f,k}\right) _{\ast }(\left[ T_{f,1},\partial T_{f,1}\right]
)
\end{eqnarray*}%
On the other hand, we also have $i_{\ast }((\phi _{f,k}^{\gamma
})_{\ast }([T_{f}^{\gamma }]))=\left( \phi _{f,k}\right) _{\ast
}(\left[ T_{f,1},\partial T_{f,1}\right] ),$ and since $i_{\ast
}:H_{3}(F/F_{k})\rightarrow H_{3}(F/F_{k},\zeta )$ is an
isomorphism, we must have $\left[ \bar{c}_{f}\right] =(\phi
_{f,k}^{\gamma })_{\ast }([T_{f}^{\gamma }]).$

Finally, notice that our choices of $\sigma \in C_{2}(\Sigma )$
and $ c_{f}\in C_{3}(\Sigma )$ certainly qualify as choices for
$\sigma \in C_{2}(F)$ and $c_{f}\in C_{3}(F),$ respectively, in
the construction of Morita's homomorphism in Section \ref{(MR)}.
Thus we have $\left( \Phi \circ \sigma _{k}\right) (f)=\Phi
(T_{f}^{\gamma },\phi _{f,k}^{\gamma })=(\phi _{f,k}^{\gamma
})_{\ast }([T_{f}^{\gamma }])=\left[ \bar{c}_{f}\right] =
\tilde{\tau}_{k}(f),$ and the theorem is proved.
\end{proof}

Now that we see that $\sigma _{k}:\mathcal{J}(k)\rightarrow \Omega
_{3}(F/F_{k})$ and Morita's homomorphism are indeed equivalent, we
can describe in a different way how $\sigma _{k}$ relates to
Johnson's homomorphism $\tau _{k}:\mathcal{J}(k)\rightarrow
H_{1}\otimes F_{k}/F_{k+1}$ . Recall the differential
\begin{equation*}
d^{2}:H_{3}\left( \frac{F}{F_{k}}\right) \rightarrow H_{1}\otimes
\frac{F_{k}}{F_{k+1}}
\end{equation*}%
discussed in Section \ref{(MR)}. Then $\tau _{k}$ factors through
$\Omega _{3}(F/F_{k})$ so that the following diagram commutes.
\begin{equation*}
\begin{diagram}
               &                      & \Omega_{3}\left( \frac{F}{F_{k}}\right) \\
               & \ruTo^{{\sigma}_{k}} & \dTo_{d^{2}\circ \Phi}                  \\
\mathcal{J}(k) & \rTo_{\tau_{k}}      & H_{1}\otimes \frac{F_{k}}{F_{k+1}}      \\
\end{diagram}
\end{equation*}

\section{A Spin Bordism Representation of the Mapping Class Group\label{(SBRMCG)}}

We introduced in Section \ref{(BRMCG)} a new representation
$\sigma _{k}: \mathcal{J}(k)\rightarrow \Omega _{3}(F/F_{k})$
which we then showed was equivalent to Morita's homomorphism
$\tilde{\tau}_{k}:\mathcal{J}(k)\rightarrow H_{3}(F/F_{k}).$
Because of the range of the latter homomorphism, it may seem
preferable to those who have a firm understanding of homology.
However, $\sigma _{k}$ has its advantages. First, it simply has a
more geometric nature to it. Second, and perhaps most importantly,
it naturally leads to an interesting question that
$\tilde{\tau}_{k}$ does not. What happens when we add more
structure to the bordism group? More specifically, what is the
result of replacing the bordism group $\Omega _{3}(F/F_{k})$ with
the spin bordism group $\Omega _{3}^{spin}(F/F_{k})$?

\subsection{A Spin Bordism Invariant of $\mathcal{J}(k)$}

Recall that a spin structure can be thought of as a trivialization
of the stable tangent bundle restricted to the 2-skeleton, and
every oriented 3-manifold has a spin structure. Since a spin
structure on a manifold induces a spin structure on its boundary,
we can define the \textit{3-dimensional spin bordism group
}$\Omega _{3}^{spin}(X)$ in exactly the same way as the oriented
bordism group $\Omega _{3}(X)$ with the additional requirement
that spin structures on spin bordant 3-manifolds must extend to a
spin structure on the 4-dimensional bordism between them. That is,
elements of $\Omega _{3}^{spin}(X)$ are equivalence classes of
triples $ \left( M,\phi ,\sigma \right) $ consisting of a closed,
spin 3-manifold $M$ with spin structure $\sigma $ and a continuous
map $\phi :M\rightarrow X.$ We say two elements $\left( M_{0},\phi
_{0},\sigma _{0}\right) $ and $\left( M_{1},\phi _{1},\sigma
_{1}\right) $ are equivalent, or \textit{spin bordant over }$X$,
if there is a triple $\left( W,\Phi ,\sigma \right) $ consisting
of a compact, spin 4-manifold $\left( W,\sigma \right) $ with
boundary $\partial (W,\sigma )=\left( M_{0},\sigma _{0}\right)
\amalg -\left( M_{1},\sigma _{1}\right) $ and a continuous map
$\Phi :W\rightarrow X$ satisfying $\Phi |_{M_{i}}=\phi _{i}.$

Further recall that the spin structures for a spin manifold $M$
are enumerated by $H^{1}(M;\mathbb{Z}_{2}).$ Thus, for example,
the number of possible spin structures on an oriented surface
$\Sigma _{g,1}$ of genus $g$ with one boundary component is
$\left\vert H^{1}(\Sigma _{g,1};\mathbb{Z}_{2})\right\vert
=\left\vert \mathbb{Z}_{2}^{2g}\right\vert =2^{2g}.$ If we fix a
spin structure on $\Sigma _{g,1}$, then we can extend it to the
product $\Sigma _{g,1}\times \left[ 0,1\right] .$ Now consider the
mapping class $f\in \mathcal{J}(k)$ for $\Sigma _{g,1}.$ For
$k\geq 2,$ $f$ acts trivially on $ H_{1}(\Sigma
_{g,1};\mathbb{Z}_{2})$ and on the set of spin structures. Thus
the spin structure on $\Sigma _{g,1}\times \left[ 0,1\right] $ can
be extended to the mapping torus $T_{f,1}.$ The number of possible
spin structures for $T_{f,1}$ is $\left\vert
H^{1}(T_{f,1};\mathbb{Z} _{2})\right\vert =\left\vert
\mathbb{Z}_{2}^{2g+1}\right\vert =2^{2g+1},$ where the extra
factor of $2$ corresponds to the extra generator $\gamma \in \pi
_{1}(T_{f,1}).$ Remember that we construct $T_{f}^{\gamma }$ from
$ T_{f,1}$ by performing a Dehn filling along $\gamma ,$ i.e.
filling the boundary $\partial T_{f,1}=\partial \Sigma
_{g,1}\times S^{1}$ with $\partial \Sigma _{g,1}\times D^{2}.$
Then, as long as we choose the spin structure for $\gamma $ which
extends over a disk, we can extend the spin structure on $T_{f,1}$
to a spin structure $\sigma $ on $T_{f}^{\gamma }.$ Again, the
number of possible spin structures for $T_{f}^{\gamma }$ is $
\left\vert H^{1}(T_{f}^{\gamma };\mathbb{Z}_{2})\right\vert
=\left\vert \mathbb{Z}_{2}^{2g}\right\vert =2^{2g},$ and these
exactly correspond to the spin structures on $\Sigma _{g,1}.$ Let
$\phi _{f,k}^{\gamma }:T_{f}^{\gamma }\rightarrow K(F/F_{k},1)$ be
as before.

\begin{theorem}
\label{(etahomom)}Let $\Sigma _{g,1}$ be a surface of genus $g$
with one boundary component and a fixed spin structure. Let
$\sigma $ denote the corresponding spin structure on
$T_{f}^{\gamma }$ for all $f\in \mathcal{J}(k)$, $k\geq 2.$ Then
there is a finite family of well-defined homomorphisms
\begin{equation*}
\eta _{\sigma,k}:\mathcal{J}(k)\rightarrow \Omega
_{3}^{spin}\left( \frac{F}{F_{k}} \right)
\end{equation*}%
defined by $\eta _{\sigma,k}(f)=\left( T_{f}^{\gamma },\phi
_{f,k}^{\gamma },\sigma \right) $.
\end{theorem}

\begin{proof}
This follows directly from the proof that $\sigma _{k}:\mathcal{J}
(k)\rightarrow \Omega _{3}(F/F_{k})$ is a well-defined
homomorphism (see Theorem \ref{(SigmaHomom)}) since the spin
structure on $T_{f}^{\gamma }\amalg T_{h}^{\gamma }$ naturally
extends over the product $\left( T_{f}^{\gamma }\amalg
T_{h}^{\gamma }\right) \times \left[ 0,1\right] $ and the spin
structure on $\Sigma _{g,1}$ naturally extends over the product
$\Sigma _{g,1}\times \left[ -\varepsilon ,\varepsilon \right]
\times \left[ -\delta ,\delta \right] .$
\end{proof}

First, we point out that if we compose this homomorphism with a
\textquotedblleft forgetful\textquotedblright\ map which ignores
the spin structure then we obtain our original homomorphism
$\sigma _{k}$. Second, recall the proof of Corollary \ref{(IO)}
where we pointed out that, by using the Atiyah-Hirzebruch spectral
sequence, one could build the $n$-dimensional bordism group
$\Omega _{n}(X,A)$ using $H_{p}(X,A;\Omega _{q})$ with $p+q=n$ as
building blocks. In the same way, the $n$-dimensional spin bordism
group $ \Omega _{n}^{spin}(X,A)$ is constructed out of
$H_{p}(X,A;\Omega _{q}^{spin}) $ with $p+q=n,$ where $\Omega
_{q}^{spin}=\Omega _{q}^{spin}(\cdot )$ is the spin bordism group
of a single point. Unlike the previous case, all but one of these
coefficient groups are nontrivial for $n=3$. In particular, since
$\Omega _{0}^{spin}\cong \mathbb{Z}$, $ \Omega _{1}^{spin}\cong
\Omega _{2}^{spin}\cong \mathbb{Z}_{2},$ and $\Omega
_{3}^{spin}\cong \left\{ e\right\} $, we have that $\Omega
_{3}^{spin}(F/F_{k})$ is built out of
\begin{equation*}
\begin{array}{ccccc}
H_{3}(F/F_{k};\Omega _{0}^{spin}) & \cong & H_{3}(F/F_{k}) & \cong & \Omega
_{3}(F/F_{k}), \\
H_{2}(F/F_{k};\Omega _{1}^{spin}) & \cong & H_{2}(F/F_{k})\otimes
\mathbb{Z}_{2} & \cong & F_{k}/F_{k+1}\otimes \mathbb{Z}_{2}, \\
H_{1}(F/F_{k};\Omega _{2}^{spin}) & \cong & H_{1}(F/F_{k})\otimes
\mathbb{Z}_{2} & \cong & \mathbb{Z}_{2}^{2g},\text{ and} \\
H_{0}(F/F_{k};\Omega _{3}^{spin}) & \cong & 0. &  &
\end{array}%
\end{equation*}%
And so at the very least we see that there is potential for $\eta
_{\sigma,k}$ to give much more information about the structure of
the group $\mathcal{J}(k).$

\subsection{A Closer Look at $\protect\eta _{\sigma,2}$\label{(CLatETA2)}}

We will now investigate the specific case when $k=2$ and see what
information $\eta _{\sigma,2}:\mathcal{J}(2)\rightarrow \Omega
_{3}^{spin}(F/F_{2})$ gives us about the Torelli group
$\mathcal{J}(2).$ We have already seen that the original Johnson
homomorphism $\tau _{2}$ factors through $\Omega
_{3}^{spin}(F/F_{2})$ (see Theorem \ref{(psi_well-def)}.) In this
section we will see that, in fact, the Birman-Craggs homomorphisms
$\left\{ \rho _{q}\right\} $ also factor through $\Omega
_{3}^{spin}(F/F_{2}).$ Therefore, this new homomorphism $\eta
_{\sigma,2}$ combines the Johnson homomorphism and Birman-Craggs
homomorphism into a single one.

Consider any mapping class $f\in \mathcal{J}(2)$ and fix a spin
structure on $\Sigma _{g,1}.$ Let $\sigma $ be the corresponding
spin structure on $ T_{f}^{\gamma }.$ Finally let $\phi
_{f}^{\gamma }=\phi _{f,2}^{\gamma }:T_{f}^{\gamma }\rightarrow
K(F/F_{2},1)$ be a continuous map which induces the canonical
epimorphism $\pi _{1}(T_{f}^{\gamma })\twoheadrightarrow
F/F_{2}\cong \mathbb{Z}^{2g}.$ Then the image under $\eta
_{\sigma,2}$ of $f$ is $ \left( T_{f}^{\gamma },\phi_{f}^{\gamma
},\sigma\right).$

The group $\left[ T_{f}^{\gamma },S^{1}\right] $ of homotopy
classes of maps $T_{f}^{\gamma }\rightarrow S^{1}$ is in
one-to-one correspondence with $\Hom(\pi _{1}(T_{f}^{\gamma
}),\mathbb{Z)}.$ In fact, there is an isomorphism $\left[
T_{f}^{\gamma },S^{1}\right] \cong H^{1}(T_{f}^{\gamma
};\mathbb{Z}).$ Let $\alpha \in H^{1}(T_{f}^{\gamma };\mathbb{Z})$
be a primitive cohomology class, then there is a continuous map
$\psi _{\alpha }:T_{f}^{\gamma }\rightarrow S^{1}$ corresponding
to $\alpha .$ There is a connected surface $S$ embedded in
$T_{f}^{\gamma }$ which represents a class in $
H_{2}(T_{f}^{\gamma })$ Poincar\'{e} dual to $\alpha ,$ and this
surface $S$ represents the same homology class in
$H_{2}(T_{f}^{\gamma })$ as $\psi _{\alpha }^{-1}(p)$ does, where
$p\in S^{1}$ is a regular value of $\psi _{\alpha }.$ (If $p\in
S^{1}$ is a regular value of $\psi _{\alpha },$ then $ \psi
_{\alpha }^{-1}(p)$ is an embedded, codimension 1 submanifold of $
T_{f}^{\gamma }$. That is, $\psi _{\alpha }^{-1}(p)$ is an
embedded surface in $T_{f}^{\gamma }$.)
\begin{figure}[h]
 \centering
 \includegraphics{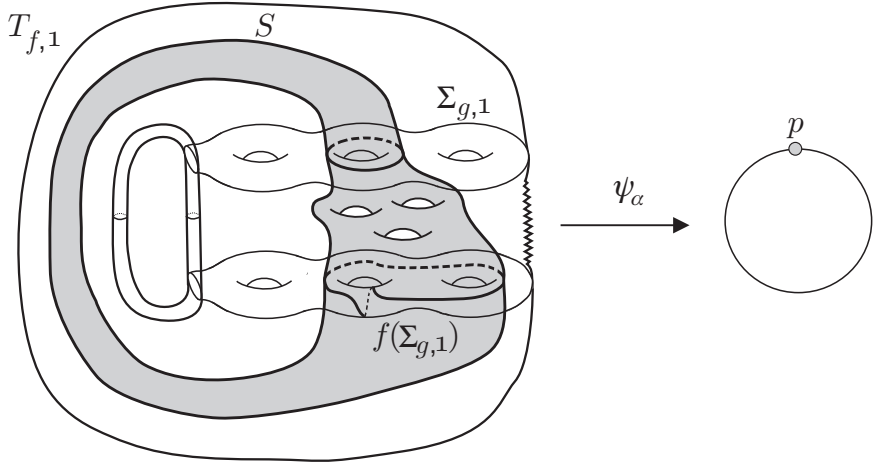}
 \caption{Embedding of $S$ into $T_{f,1}\hookrightarrow T_{f}^{\protect\gamma }$ and
the map $\protect\psi _{\protect\alpha }$}
 \label{fig emb surface}
\end{figure}

Let $\pi _{\alpha }:K(F/F_{2},1)\rightarrow S^{1}$ be a continuous
map such that $\psi _{\alpha }$ is homotopic to $\pi _{\alpha
}\circ \phi _{f}^{\gamma }$, and let $\left( \pi _{\alpha }\right)
_{\ast }:\Omega _{3}^{spin}(F/F_{2})\rightarrow \Omega
_{3}^{spin}(S^{1})$ denote the induced bordism homomorphism. Then
we can define a homomorphism
\begin{equation*}
\omega _{\sigma,\alpha }=\left( \pi _{\alpha }\right) _{\ast
}\circ \eta _{\sigma,2}: \mathcal{J}(2)\rightarrow \Omega
_{3}^{spin}(S^{1})
\end{equation*}%
by sending $f\in \mathcal{J}(2)$ to the bordism class $\left(
T_{f}^{\gamma },\psi _{\alpha },\sigma \right) \in \Omega
_{3}^{spin}(S^{1}).$ Again, using the Atiyah-Hirzebruch spectral
sequence, we see that $\Omega _{3}^{spin}(S^{1})\cong \Omega
_{2}^{spin}\cong \mathbb{Z}_{2}.$ The specific isomorphism is
given by $\left( M,\phi ,\sigma \right) \mapsto \left( \phi
^{-1}(p),\sigma |_{\phi ^{-1}(p)}\right) $, where $p\in S^{1}$ is
a regular value of $\phi .$ We can see by this isomorphism that
the spin structure on $T_{f}^{\gamma }$ restricts to a spin
structure on $S=\psi _{\alpha }^{-1}(p).$

\begin{theorem}
\label{(MyBCH)}The fixed spin structure $\sigma$ on $\Sigma
_{g,1}$ has a canonically associated quadratic form
$q:H_{1}(\Sigma _{g,1};\mathbb{Z}_{2})\rightarrow \mathbb{Z}_{2}$.
If $\arf(\Sigma _{g,1},q)=0$, there is a primitive cohomology
class $\alpha \in H^{1}(T_{f}^{\gamma };\mathbb{Z})$ such that the
homomorphism $\omega _{\sigma,\alpha }:\mathcal{J}(2)\rightarrow
\Omega _{3}^{spin}(S^{1})$ is equivalent to the Birman-Craggs
homomorphism $ \rho _{q}:\mathcal{J}(2)\rightarrow
\mathbb{Z}_{2}$.
\end{theorem}

We note that the hypothesis $\arf(\Sigma _{g,1},q)=0$ is necessary
for the Birman-Craggs homomorphism $\rho
_{q}:\mathcal{J}(2)\rightarrow \mathbb{Z }_{2}$ to be defined. See
Section \ref{(BCH)} for details.

We have a surface $S=\psi _{\alpha }^{-1}(p)$ embedded in
$T_{f}^{\gamma }.$ To determine whether the image of $f$ under the
homomorphism $\omega _{\sigma,\alpha }:\mathcal{J}(2)\rightarrow
\Omega _{3}^{spin}(S^{1})$ is trivial or not, we simply need to
determine $\left( S,\sigma |_{S}\right) \in \Omega _{2}^{spin}$.
However, this is just the well-known Arf invariant of $S$ with
respect to $\sigma |_{S}.$ We defined the Arf invariant
$\arf(\Sigma ,q)$ for a closed surface $\Sigma $ and
$\mathbb{Z}_{2}$-quadratic form $q:H_{1}(\Sigma
;\mathbb{Z}_{2})\rightarrow \mathbb{Z}_{2}$ in Section \ref
{(BCH)}. For the spin structure $\sigma |_{S}$ on $S$ let
$q_{\sigma }:H_{1}(S;\mathbb{Z}_{2})\rightarrow \mathbb{Z}_{2}$ be
the corresponding $ \mathbb{Z}_{2}$-quadratic form. Namely,
$q_{\sigma }$ is defined to be the quadratic form given by
$q_{\sigma }(x)=0$ if $\sigma |_{x}$ is the spin structure that
extends over a disk and $q_{\sigma }(x)=1$ if $\sigma |_{x}$ does
not extend over a disk. It is the work of Johnson in \cite{[J1]}
that tells us this quadratic form is equivalent to the quadratic
form discussed in Section \ref{(BCH)}. Then we have
\begin{equation*}
\arf(S,q_{\sigma })=\arf(S,\sigma |_{S})=\left( S,\sigma
|_{S}\right) \in \Omega _{2}^{spin}.
\end{equation*}

We will also need a more general definition of the Arf invariant
which includes surfaces with boundary. The definition is the same
except for a small change to the $\mathbb{Z}_{2}$-quadratic form
$q.$ In particular we have a $\mathbb{Z}_{2}$-quadratic form
\begin{equation*}
q:\frac{H_{1}(\Sigma ;\mathbb{Z}_{2})}{i_{\ast }(H_{1}(\partial
\Sigma ; \mathbb{Z}_{2}))}\rightarrow \mathbb{Z}_{2}
\end{equation*}%
where $i_{\ast }$ is induced by inclusion $i:\partial \Sigma
\rightarrow \Sigma .$ Then for a symplectic basis $\left\{
x_{i},y_{i}\right\} $ of the quotient $H_{1}(\Sigma
;\mathbb{Z}_{2})/i_{\ast }(H_{1}(\partial \Sigma ;
\mathbb{Z}_{2}))$, the \textit{Arf invariant }of $\Sigma $ with
respect to $ q $ is defined to be
\begin{equation*}
\arf(\Sigma ,q)=\sum_{i=1}^{g}q(x_{i})q(y_{i})\text{ (mod 2).}
\end{equation*}%
Notice that if the surface $\Sigma $ happens to be embedded in
$S^{3}$ then this definition is the same as the definition of the
\textit{Arf invariant }$ \arf(L)$ \textit{of an oriented link} $L$
in $S^{3}$ with components $ \left\{ L_{i}\right\} $ and
satisfying the property that the linking number is $\lk\left(
L_{i},L-L_{i}\right) \equiv 0$ modulo 2. The surface $ \Sigma $
would be a Seifert surface for the link, and $q$ would be the mod
2 Seifert self-linking form on $H_{1}(\Sigma
;\mathbb{Z}_{2})/i_{\ast }(H_{1}(\partial \Sigma
;\mathbb{Z}_{2})),$ where the self-linking is computed with
respect to a push-off in a direction normal to the surface. See
the W. Lickorish text \cite{[Li]} for more details.

Now consider the surface $S=\psi _{\alpha }^{-1}(p)$ embedded in
$T_{f}^{\gamma },$ and suppose that $S$ has genus $k$. There
exists a symplectic basis $\left\{ x_{i},y_{i}\right\} ,$ $1\leq
i\leq k,$ of $ H_{1}(S;\mathbb{Z}_{2})$ such that $x_{k}$ is
homologous to the homology class $\left[ \gamma \right] $
corresponding to $\gamma $ in $T_{f}^{\gamma } $ and $y_{k}$ is
homologous to the homology class of $\beta =S\cap \Sigma
_{g,1}\subset T_{f}^{\gamma }.$ But $\gamma $ was required to have
the spin structure that extends over a disk (so the spin structure
on $T_{f,1}$ may be extended to a spin structure on $T_{f}^{\gamma
}$.) Thus $q_{\sigma }(x_{k})=q_{\sigma }(\left[ \gamma \right]
)=0,$ and
\begin{equation*}
\arf(S,q_{\sigma })=\sum_{i=1}^{k}q_{\sigma }(x_{i})q_{\sigma
}(y_{i})=\sum_{i=1}^{k-1}q_{\sigma }(x_{i})q_{\sigma }(y_{i}).
\end{equation*}%

If we cut $S$ open along a simple closed curve parallel to $\beta
=S\cap \Sigma _{g,1}$ then the result deformation retracts to a
surface $S^{\prime } $ with boundary $\partial S^{\prime }=\beta
\amalg -f(\beta )$ and such that $ H_{1}(S^{\prime
};\mathbb{Z}_{2})$ has symplectic basis $\left\{
x_{i},y_{i}\right\}$, $1\leq i\leq k-1$. (See Figure \ref{fig
cutopen}.) If we let $q_{\sigma }^{\prime }:H_{1}(S^{\prime
};\mathbb{Z}_{2})/i_{\ast }(H_{1}(\partial S^{\prime };
\mathbb{Z}_{2}))\rightarrow \mathbb{Z}_{2}$ be the induced
$\mathbb{Z}_{2}$ -quadratic form, then
\begin{equation*}
\arf(S,q_{\sigma })=\sum_{i=1}^{k-1}q_{\sigma
}(x_{i})q_{\sigma}(y_{i})=\sum_{i=1}^{k-1}q_{\sigma }^{\prime
}(x_{i})q_{\sigma }^{\prime }(y_{i})=\arf(S^{\prime },q_{\sigma
}^{\prime }).
\end{equation*}%
According to Johnson in \cite{[J1]}, the quadratic form $q$ (in
the statement of Theorem \ref{(MyBCH)}) corresponds to a Heegaard
embedding of $ \Sigma _{g,1}$ into $S^{3}$. Thus we get an induced
embedding of $\Sigma _{g,1}\times \left[ 0,1\right] $, and thus of
$S^{\prime }$, into $S^{3}$, and the quadratic form $q_{\sigma
}^{\prime }$ is precisely the same as the mod 2 Seifert
self-linking form. Thus we see that to calculate $\arf
(S,q_{\sigma })$, we really only need to calculate the Arf
invariant of the link $\left\{ \beta ,f(\beta )\right\} $ with
Seifert surface $S^{\prime }.$
\begin{figure}[h]
 \centering
 \includegraphics{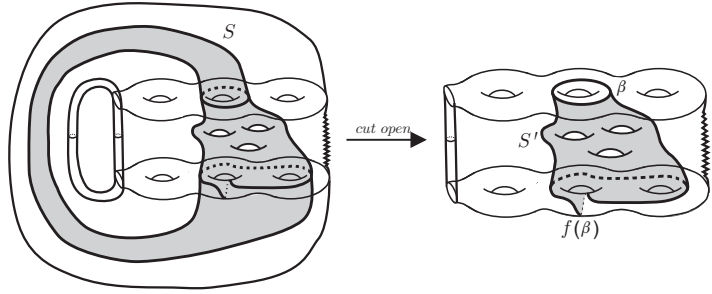}
 \caption{$S$ cut open along $\protect\beta $ to obtain $S^{\prime }$}
 \label{fig cutopen}
\end{figure}

Since there is an isomorphism $H^{1}(T_{f}^{\gamma
};\mathbb{Z})\cong H^{1}(\Sigma _{g,1};\mathbb{Z)},$ $\alpha \in
H^{1}(T_{f}^{\gamma };\mathbb{Z})$ has a corresponding class in
$H^{1}(\Sigma _{g,1};\mathbb{Z)}$ which we will also call $\alpha
.$ The homology class of $\beta =S\cap \Sigma _{g,1}$ in
$H_{1}(T_{f}^{\gamma })$ also has a corresponding class $\left[
\beta \right]$ in $H_{1}(\Sigma _{g,1}).$ Since the homology class
of $S$ is Poincar\'{e} dual to $\alpha \in H^{1}(T_{f}^{\gamma };
\mathbb{Z}),$ $\left[ \beta \right] \in H_{1}(\Sigma _{g,1})$ must
be Poincar\'{e} dual to $\alpha \in H^{1}(\Sigma
_{g,1};\mathbb{Z)}$.

\begin{proof}[Proof of Theorem \protect\ref{(MyBCH)}]
We have a fixed spin structure $\sigma$ on $\Sigma _{g,1}.$ Let
$q$ be the associated $\mathbb{Z}_{2}$-quadratic form. Recall from
Section \ref{(BCH)} that the hypothesis $\arf(\Sigma _{g,1},q)=0$
was necessary for the Birman-Craggs homomorphism $\rho
_{q}:\mathcal{J}(2)\rightarrow \mathbb{Z} _{2}$ to be defined. We
have already seen that the spin structure on $\Sigma _{g,1}$
induces a spin structure on $T_{f}^{\gamma }$ which we will also
denote by $\sigma $ and which in turn induces a spin structure
$\sigma |_{S}$ on the surface $S$ defined above. To prove the
theorem, we need to find a primitive cohomology class $ \alpha \in
H^{1}(T_{f}^{\gamma };\mathbb{Z})$ such that $\omega
_{\sigma,\alpha }(f)=\rho _{q}(f).$ To accomplish this we need to
find a surface $S$ that represents a homology class Poincar\'{e}
dual to $ \alpha$ and such that $\arf(S,q_{\sigma
})=\arf(S^{\prime },q_{\sigma }^{\prime })=\rho _{q}(f).$ To do
so, we will construct a simple closed curve $\beta $ on $\Sigma
_{g,1}$ and calculate the Arf invariant $\arf(\beta ,f(\beta ))$
with Seifert surface $S^{\prime }$ in $\Sigma _{g,1}\times \left[
0,1\right] \hookrightarrow S^{3}.$

Recall that for genus $g=2$ surfaces, the Torelli group
$\mathcal{J}(2)$ is generated by the collection of all Dehn twists
about bounding simple closed curves, and for genus $g\geq 3$,
$\mathcal{J}(2)$ is generated by the collection of all Dehn twists
about genus 1 cobounding pairs of simple closed curves, i.e. pairs
of non-bounding, disjoint, homologous simple closed curves that
together bound a genus 1 subsurface. Thus it is sufficient to
prove the claim for such elements of $\mathcal{J}(2).$

First assume that $g=2$ and $C$ is a genus 1 bounding simple
closed curve on $\Sigma _{2\text{,}1}$. Let $f$ be a Dehn twist
about $C.$ Then $C$ splits $ \Sigma _{2,1}$ into two genus 1
surfaces $\Sigma _{a}$ and $\Sigma _{b}.$ Let $\left\{
x_{a},y_{a}\right\} $ and $\left\{ x_{b},y_{b}\right\} $ be
symplectic bases of $H_{1}(\Sigma _{a})/i_{\ast }(H_{1}(\partial
\Sigma _{a}))$ and $H_{1}(\Sigma _{b})/i_{\ast }(H_{1}(\partial
\Sigma _{b})),$ respectively. Then we have two cases:

\begin{enumerate}
\item[(i)] $\rho _{q}(f)=\arf(\Sigma _{a},q|_{\Sigma
_{a}})=\arf(\Sigma _{b},q|_{\Sigma _{b}})=1$

$\iff $ $q(x_{a})=q(y_{a})=q(x_{b})=q(y_{b})=1,$ or

\item[(ii)] $\rho _{q}(f)=\arf(\Sigma _{a},q|_{\Sigma
_{a}})=\arf(\Sigma _{b},q|_{\Sigma _{b}})=0$

$\iff $ at least one of $\left\{q(x_{a}),q(y_{a})\right\} $ and
one of $\left\{q(x_{b}),q(y_{b})\right\} $ are $0.$
\end{enumerate}

Without loss of generality, let us assume in case (ii) that
$q(x_{a})=q(x_{b})=0.$ Then in either case we have $\rho
_{q}(f)=q(x_{a}).$ Let $\beta $ be a simple closed curve on
$\Sigma _{2,1}\hookrightarrow T_{f}^{\gamma }$ which intersects
$C$ exactly twice and such that $\left[ \beta \right] \in
H_{1}(\Sigma _{2\text{,}1})$ is homologous to $x_{a}+x_{b}$. Then
we also have the simple closed curve $f(\beta )$ on $f(\Sigma
_{2,1})\hookrightarrow T_{f}^{\gamma }.$ Near $C$ the picture will
always be as in Figure \ref{fig twist1}, and we choose $ S^{\prime
}$ to be this particular surface pictured in Figure \ref{fig
twist1} with boundary $\partial S^{\prime }=\beta \amalg -f(\beta
).$
\begin{figure}[h]
 \centering
 \includegraphics{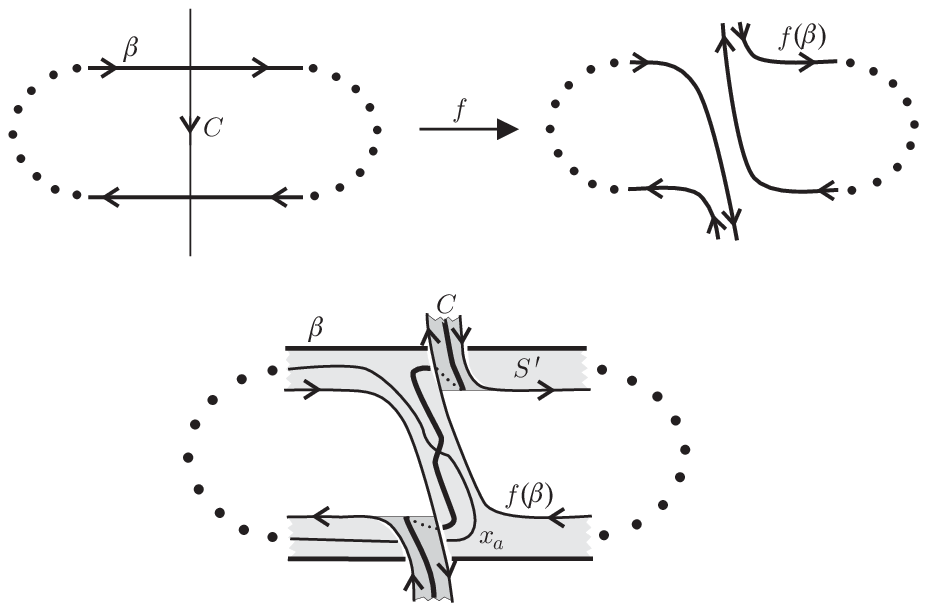}
 \caption{Surface $S^{\prime }$ in $T_{f}^{\protect \gamma }$ with
 boundary $\protect\beta \amalg -f\left(\protect\beta \right) $ (for $g=2$)}
 \label{fig twist1}
\end{figure}
This surface $S^{\prime }$ has spin structure $\sigma |_{S^{\prime
}}$ and a corresponding quadratic form
\begin{equation*}
q_{\sigma }^{\prime }:H_{1}(S^{\prime };\mathbb{Z}_{2})/i_{\ast
}(H_{1}(\partial S^{\prime };\mathbb{Z}_{2}))\rightarrow \mathbb{Z}_{2}
\end{equation*}%
given by the mod 2 self-linking form. Notice that $\left\{
x_{a},\left[ C \right] \right\} $ is a symplectic basis for the
quotient $H_{1}(S^{\prime }; \mathbb{Z}_{2})/i_{\ast
}(H_{1}(\partial S^{\prime };\mathbb{Z}_{2})).$ Then we have
\begin{equation*}
\omega _{\sigma,\alpha }(f)=\arf(S^{\prime },q_{\sigma }^{\prime
})\overset{ def.}{\equiv }\arf(\beta ,f(\beta ))=q_{\sigma
}^{\prime }(x_{a})q_{\sigma }^{\prime }(\left[ C\right] ).
\end{equation*}%
Note that, while $C$ is a product of commutators on $\Sigma
_{2,1},$ it is not a product of commutators on $S^{\prime }.$ But
it is easy to see from Figure \ref{fig twist1} that $q_{\sigma
}^{\prime }(\left[ C\right] )=\lk \left( C,C^{+}\right) \equiv 1$
modulo 2. It is also clear that $q_{\sigma }^{\prime
}(x_{a})=q(x_{a}).$ Thus
\begin{equation*}
\omega _{\sigma,\alpha }(f)=q_{\sigma}^{\prime }(x_{a})q_{\sigma
}^{\prime }(\left[ C\right] )=q(x_{a})=\arf(\Sigma _{a},q|_{\Sigma
_{a}})=\rho _{q}(f).
\end{equation*}

Now assume that $g\geq 3$ and $C_{1}$ and $C_{2}$ are genus 1
cobounding pairs of simple closed curves on $\Sigma _{g,1}.$ Let
$f$ be a composition of Dehn twists about $C_{1}$ and $C_{2}.$
Then $C_{1}$ and $C_{2}$ cobound a genus 1 subsurface $\Sigma
^{\prime }.$ Let $\left\{ x,y\right\} $ be a symplectic basis of
$H_{1}(\Sigma ^{\prime })/i_{\ast }(H_{1}(\partial \Sigma ^{\prime
})).$ There are two cases:

\begin{enumerate}
\item[(1)] $q(C_{1})=q(C_{2})=1$ and

\item[(2)] $q(C_{1})=q(C_{2})=0.$
\end{enumerate}

For case (1), we simply let $\beta $ be a simple closed curve on
$\Sigma _{g,1}\hookrightarrow T_{f}^{\gamma }$ which does not
intersect $C_{1}$ or $ C_{2}$. Then $f$ will not affect $\beta ,$
and we can choose $S^{\prime }$ to be a straight cylinder between
$\beta $ and $f(\beta )$ so that $ H_{1}(S^{\prime
};\mathbb{Z}_{2})/i_{\ast }(H_{1}(\partial S^{\prime };
\mathbb{Z}_{2}))$ is trivial. Thus $\omega _{\sigma,\alpha
}(f)=\arf (S^{\prime },q_{\sigma }^{\prime })=0.$ We also know
from the end of Section \ref{(BCH)} that in this case $\rho
_{q}(f)=0$.

For case (2), we have two subcases:

\begin{enumerate}
\item[(i)] $\rho _{q}(f)=\arf(\Sigma^{\prime },q|_{\Sigma ^{\prime
}})=1\iff $ $q(x)=q(y)=1,$ or

\item[(ii)] $\rho _{q}(f)=\arf(\Sigma ^{\prime
},q|_{\Sigma^{\prime }})=0\iff $ at least one of $\left\{
q(x),q(y)\right\} $ is $0.$
\end{enumerate}

Again without loss of generality, let us assume in case (ii) that
$q(x)=0.$ In both cases let $\beta $ be a simple closed curve on
$\Sigma _{g,1}\hookrightarrow T_{f}^{\gamma }$ which intersects
each of $C_{1}$ and $ C_{2}$ exactly once and such that $\left[
\beta \right] \in H_{1}(\Sigma _{g,1})$ is homologous to
$x+x^{\prime },$ where $x^{\prime }$ is any nontrivial homology
class in $H_{1}(\Sigma _{g,1}-\Sigma ^{\prime })$. Then we also
have the simple closed curve $f(\beta )$ on $f(\Sigma
_{g,1})\hookrightarrow T_{f}^{\gamma }.$ Near $C_{1}$ and $C_{2}$
the picture will always be as in Figure \ref{fig twist2}, and we
choose $S^{\prime }$ to be this particular surface pictured in
Figure \ref{fig twist2} with boundary $\partial S^{\prime }=\beta
\amalg -f(\beta ).$
\begin{figure}[h]
 \centering
 \includegraphics{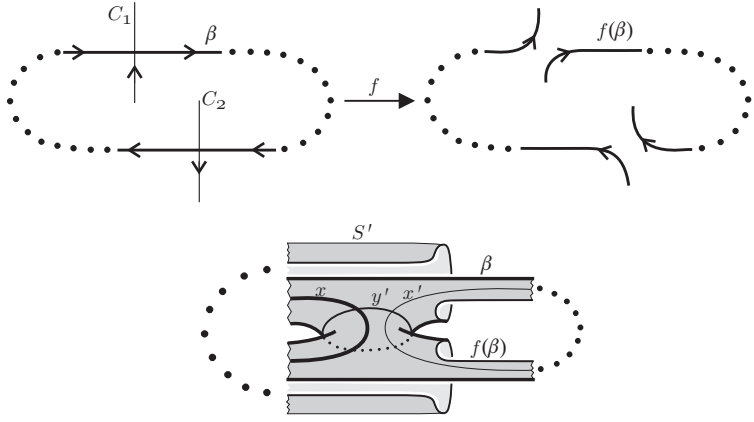}
 \caption{Surface $S^{\prime }$ in $T_{f}^{\protect\gamma }$ with
 boundary $\protect\beta \amalg -f\left(\protect\beta \right) $
 (for $g\geq 3$)}
 \label{fig twist2}
\end{figure}
Again this surface $S^{\prime }$ has spin structure $\sigma
|_{S^{\prime }}$ and a corresponding quadratic form $q_{\sigma
}^{\prime }:H_{1}(S^{\prime };\mathbb{Z}_{2})/i_{\ast
}(H_{1}(\partial S^{\prime }; \mathbb{Z}_{2}))\rightarrow
\mathbb{Z}_{2}$ given by the mod 2 self-linking form. Let
$y^{\prime }$ be any homology class such that $\left\{ x,y^{\prime
}\right\} $ is a symplectic basis for $H_{1}(S^{\prime
};\mathbb{Z}_{2})/i_{\ast }(H_{1}(\partial S^{\prime
};\mathbb{Z}_{2})).$ Then we have
\begin{equation*}
\omega _{\sigma,\alpha }(f)=\arf(S^{\prime },q_{\sigma }^{\prime
})\overset{def.}{\equiv }\arf(\beta ,f(\beta ))=q_{\sigma
}^{\prime }(x)q_{\sigma }^{\prime }(y^{\prime }).
\end{equation*}%
Notice that $\left\{ x,y^{\prime }\right\} $ is also a basis for
$H_{1}(\Sigma ^{\prime };\mathbb{Z}_{2})/i_{\ast }(H_{1}(\partial
\Sigma ^{\prime };\mathbb{Z}_{2}))$ and that $q_{\sigma }^{\prime
}(x)=q(x)$ and $ q_{\sigma }^{\prime }(y^{\prime })=q(y^{\prime
}).$ Thus we see that
\begin{equation*}
\omega _{\sigma,\alpha }(f)=q_{\sigma }^{\prime }(x)q_{\sigma
}^{\prime }(y^{\prime })=q(x)q(y^{\prime })=\arf(\Sigma ^{\prime
},q|_{\Sigma ^{\prime }})=\rho _{q}(f).
\end{equation*}%
This completes the proof of Theorem \ref{(MyBCH)}.
\end{proof}

As a result of this theorem and Theorem \ref{(psi_well-def)}, we
see that $ \eta _{\sigma,2}$ contains the necessary information
for determining both the Johnson homomorphism $\tau _{2}$ and the
Birman-Craggs homomorphism $\rho _{q}.$ Recall from Section
\ref{(ATG)} that the abelianization
$H_{1}(\mathcal{J}(2);\mathbb{Z)}\cong \mathcal{J}(2)/\left[
\mathcal{J}(2),\mathcal{J}(2)\right] $ of the Torelli group is
completely determined by $\tau _{2}$ and $\rho _{q}$ (over all
possible $q$) since the commutator subgroup of the Torelli group
is given by the kernels of these homomorphisms. Namely, we have
$\left[ \mathcal{J}(2),\mathcal{J}(2)\right] =\mathcal{C}\cap
\mathcal{J}(3),$ where $\mathcal{C}=\bigcap_{q}\ker \rho _{q}.$
Suppose we take a mapping class $f\in \ker \eta _{\sigma,2}.$
Certainly it is true that $f\in \ker \rho_{q} \cap \mathcal{J}(3)$
since $\tau _{2}$ and $\rho _{q} $ factor through $\eta
_{\sigma,2}.$ Moreover, if
\begin{equation*}
\mathcal{D}=\bigcap\limits_{\sigma }\ker \eta _{\sigma,2}
\end{equation*}%
is the common kernel over all possible spin structures, then
$\mathcal{D}\subset \mathcal{C}\cap \mathcal{J}(3).$ Of course it
would be nice to know if the converse is also true.

\begin{problem}
What is $\ker \eta _{\sigma,2}$? Is $\ker \eta _{\sigma,2}=\ker
\rho_{q} \cap \mathcal{J}(3)$?
\end{problem}
\begin{problem}
Is it true that $\mathcal{D}=\mathcal{C\cap J}(3)=\left[
\mathcal{J}(2),\mathcal{J}(2)\right] $?
\end{problem}

\subsection{Analysis of $\protect\eta _{\sigma,k}$}

In this section we shift our focus to the general homomorphism
$\eta _{\sigma,k}:\mathcal{J}(k)\rightarrow \Omega
_{3}^{spin}(F/F_{k})$ for arbitrary values of $k.$ We already know
that $\ker \eta _{\sigma,k}\subset
\mathcal{J}\left(2k-1\right)=\ker \sigma _{k}$ since the oriented
bordism homomorphism $\sigma _{k}:\mathcal{J}(k)\rightarrow \Omega
_{3}(F/F_{k})$ factors through $\Omega _{3}^{spin}(F/F_{k}).$
However, the additional structure on the bordism given by the spin
structures should refine the kernel of $\eta _{\sigma,k}. $

\begin{problem}
What is the kernel of $\eta _{\sigma,k}$?
\end{problem}

\begin{problem}
Does $\eta _{\sigma,k}$ give a faithful representation of the
abelianization of $ \mathcal{J}(k)$? In other words, is
$\im\eta_{\sigma,k}\cong\mathcal{J}(k)/\left[\mathcal{J}(k),\mathcal{J}(k)\right]$?
\end{problem}

A sufficient condition for $f\in \ker \eta _{\sigma,k}$ is given
in the following theorem, but it is most likely not necessary.
Consider the entire collection $\left\{ \omega _{\sigma,\alpha
}\right\} $ of the homomorphisms $\omega _{\sigma,\alpha }:
\mathcal{J}(2)\rightarrow \Omega _{3}^{spin}(S^{1})$ defined in
Section \ref {(CLatETA2)}, and let
\begin{equation*}
\mathcal{B=}\bigcap\limits_{\alpha }\ker \omega _{\sigma,\alpha }
\end{equation*}%
be the common kernel of all $\omega_{\sigma,\alpha }$ for all
$\alpha \in H^{1}(T_{f}^{\gamma };\mathbb{Z}).$

\begin{theorem}
\label{(kerEtaK)}If $f\in \mathcal{B\cap J}(2k+1),$ then $f\in
\ker \eta _{\sigma,k}.$
\end{theorem}

Note that the hypothesis requires $f\in
\mathcal{J}\left(2k+1\right),$ not just $f\in
\mathcal{J}\left(2k-1\right).$ The purpose of this will be
revealed in the proof of the theorem, but it is probably not
necessary. However, as stated above, it is certainly necessary
that $f\in \mathcal{J}\left(2k-1\right).$

Before we give the proof of this theorem, let us first set up some
necessary notation. For a more complete discussion, we refer the
reader to Whitehead's book \cite{[Wh]}. We will be using the
Atiyah-Hirzebruch spectral sequence. In particular, let
\begin{equation}
J_{p,q}^{m}=\image\left( \left( i_{p,q}\right) _{\ast
}:\tilde{\Omega}_{p+q}^{spin}\left( \left( \frac{F}{F_{m}}\right)
^{\left( p\right) }\right) \rightarrow
\tilde{\Omega}_{p+q}^{spin}\left( \frac{F}{F_{m}}\right) \right)
\text{.}  \tag{$\star $}  \label{(Jpq)}
\end{equation}%
Here $(F/F_{m})^{(p)}$ denotes the $p$-skeleton of $K(F/F_{m},1)$,
$\left( i_{p,q}\right) _{\ast }$ is induced by the inclusion map $
(F/F_{m})^{(p)}\hookrightarrow K(F/F_{m},1)$, and $\tilde{\Omega}
_{n}^{spin}(F/F_{m})$ denotes the \textit{reduced} spin bordism
group defined by
\begin{equation*}
\Omega _{n}^{spin}\left( \frac{F}{F_{m}}\right) \cong \Omega
_{n}^{spin}\oplus \tilde{\Omega}_{n}^{spin}\left(
\frac{F}{F_{m}}\right).
\end{equation*}%
Note that if $\left( M,\phi ,\sigma \right) \in J_{p,q}^{m}$ then
for $l\leq m$ the triple $\left( M,\pi _{m,l}\circ \phi ,\sigma
\right) $ is in $ J_{p,q}^{l}$, where $\pi
_{m,l}:K(F/F_{m},1)\rightarrow K(F/F_{l},1)$ is the projection
map. Let
\begin{equation*}
E_{p,q}^{2}\cong \tilde{H}_{p}(F/F_{m};\Omega _{q}^{spin}),
\end{equation*}%
and the boundary operator is
\begin{equation*}
d_{p,q}^{2}:E_{p,q}^{2} \rightarrow E_{p-2,q+1}^{2}.
\end{equation*}%
The groups $E_{p,q}^{2}$ may be thought of as the building blocks
for $\tilde{\Omega}_{n}^{spin}(F/F_{m})$ with $p+q=n$. In
actuality, the building blocks are the groups $E_{p,q}^{\infty
}=\lim E_{p,q}^{r}$, where for $r\geq 3$
\begin{equation*}
E_{p,q}^{r}=\frac{\ker d_{p,q}^{r-1}}{\im
d_{p+r-1,q-r+2}^{r-1}}\text{ and
}d_{p,q}^{r-1}:E_{p,q}^{r-1}\rightarrow
E_{p-r+1,q+r-2}^{r-1}\text{.}
\end{equation*}%
We also have an isomorphism
\begin{equation}
E_{p,q}^{\infty }\cong J_{p,q}^{m}/J_{p-1,q+1}^{m}.  \tag{$\star
\star $} \label{(Epq)}
\end{equation}

Since $\Omega _{3}^{spin}=0$, we then have $\Omega
_{3}^{spin}(F/F_{m})\cong \tilde{\Omega}_{3}^{spin}(F/F_{m})$ and
\begin{equation*}
\Omega _{3}^{spin}(F/F_{m})=J_{3,0}^{m}\supseteq
J_{2,1}^{m}\supseteq J_{1,2}^{m}\supseteq J_{0,3}^{m}=0\text{.}
\end{equation*}%
Then one can show that the relevant $E_{p,q}^{\infty }$ are as
follows.
\begin{eqnarray*}
E_{3,0}^{\infty } &=&E_{3,0}^{3}=\ker d_{3,0}^{2}\subset \tilde{H}
_{3}(F/F_{m})\cong H_{3}(F/F_{m}) \\
E_{2,1}^{\infty } &=&E_{2,1}^{3}=\coker d_{4,0}^{2}\cong
H_{2}(F/F_{m};\Omega _{1}^{spin})/\im d_{4,0}^{2} \\
E_{1,2}^{\infty } &=&E_{1,2}^{4}=\coker d_{4,0}^{3}\cong
H_{1}(F/F_{m};\Omega _{2}^{spin})/\im d_{3,1}^{2} \\
E_{0,3}^{\infty } &=&0
\end{eqnarray*}%
We can now begin our proof of Theorem \ref{(kerEtaK)}.

\begin{proof}[Proof of Theorem \protect\ref{(kerEtaK)}]
We assume that $f\in \mathcal{B\cap J}(2k+1)$, and we want to show
that $ \left( T_{f}^{\gamma },\phi _{f,k}^{\gamma },\sigma \right)
=0$ in $\Omega _{3}^{spin}(F/F_{k})$. Since $\Omega
_{3}^{spin}(F/F_{m})=J_{3,0}^{m}$, we may perturb any $\left(
M,\phi ,\sigma \right) \in \Omega _{3}^{spin}(F/F_{m})$ to ensure
that $\phi (M)$ is contained in the 3-skeleton $(F/F_{m})^{(3)}$
of $K(F/F_{m},1)$. That is, by the definition of $J_{3,0}^{m}$
given in (\ref{(Jpq)}) we can choose $\phi ^{\prime }$ homotopic
to $\phi $ so that $\left( M,\phi ^{\prime },\sigma \right) \in
\Omega _{3}^{spin}((F/F_{m})^{(3)})$ and $\left( i_{3,0}\right)
_{\ast }(M,\phi ^{\prime },\sigma )=\left( M,\phi ,\sigma \right)
$ in $\Omega _{3}^{spin}(F/F_{m})$.

Since $f\in \mathcal{J}\left(2k+1\right)\subset
\mathcal{J}\left(k+1\right)$, Lemma \ref{(phi)} says
$\phi_{f,k+1}^{\gamma}$ exists. We start with $\left(
T_{f}^{\gamma },\phi _{f,k+1}^{\gamma },\sigma \right) \in
J_{3,0}^{k+1}$. Then Theorem \ref {(nullbordant)} says that the
pair $\left( T_{f}^{\gamma },\phi _{f,k+1}^{\gamma }\right) =0$ in
$\Omega _{3}(F/F_{k+1})$. Thus ($\phi _{f,k+1}^{\gamma })_{\ast
}([T_{f}^{\gamma }])=0$ in $H_{3}(F/F_{k+1})\cong E_{3,0}^{\infty
}$, and we therefore know from (\ref{(Epq)}) that $\left(
T_{f}^{\gamma },\phi _{f,k+1}^{\gamma },\sigma \right) $ must be
in $ J_{2,1}^{k+1}$. Thus by (\ref{(Jpq)}) there exists a triple
$\left( M^{\prime },\phi ^{\prime },\sigma ^{\prime }\right) \in
\Omega _{3}^{spin}((F/F_{k+1})^{(2)})$ such that $\left(
i_{2,1}\right) _{\ast }(M^{\prime },\phi ^{\prime },\sigma
^{\prime })=\left( T_{f}^{\gamma },\phi _{f,k+1}^{\gamma },\sigma
\right) $ in $\Omega _{3}^{spin}(F/F_{k+1})$ as indicated in the
following diagram.
\begin{equation*}
\begin{diagram}
                                                          & &          & \Omega_{3}^{spin}( (F/F & _{k+1})^{(2)}) \ \ \ \ \ \ \ \                   \\
                                                          & &          & \dTo                    & ( M^{\prime },\phi^{\prime },\sigma ^{\prime })  \\
\Omega_{3}^{spin}( F/F_{2k+1})                            & & \rTo     & \Omega_{3}^{spin}( F/F  & _{k+1}) \ \ \ \  \dMapsto                        \\
 \ \ \ (T_{f}^{\gamma },\phi _{f,2k+1}^{\gamma },\sigma ) & & \rMapsto &                         & (T_{f}^{\gamma },\phi_{f,k+1}^{\gamma },\sigma ) \\
\end{diagram}\bigskip
\end{equation*}

\begin{lemma}
\label{(zero)}The homomorphism $\left( \pi _{k+1,k}\right) _{\ast
}:J_{2,1}^{k+1}/J_{1,2}^{k+1}\longrightarrow
J_{2,1}^{k}/J_{1,2}^{k}$ is trivial.
\end{lemma}

\begin{proof}
By (\ref{(Epq)}) we have $J_{2,1}^{k+1}/J_{1,2}^{k+1}\cong
E_{2,1}^{\infty }\cong H_{2}(F/F_{k+1};\Omega _{1}^{spin})/\im
d_{4,0}^{2}$. Similarly, we have $J_{2,1}^{k}/J_{1,2}^{k}\cong
H_{2}(F/F_{k};\Omega _{1}^{spin})/ \im d_{4,0}^{2}$. So this
homomorphism is equivalent to
\begin{equation*}
\frac{H_{2}(F/F_{k+1};\Omega _{1}^{spin})}{\im
d_{4,0}^{2}}\rightarrow \frac{H_{2}(F/F_{k};\Omega
_{1}^{spin})}{\im d_{4,0}^{2}}.
\end{equation*}%
In the proof of Corollary \ref{(TuraevCor)} we showed that
$H_{2}(F/F_{k+1})\rightarrow H_{2}(F/F_{k})$ is the zero map. Thus
$H_{2}(F/F_{k+1};\Omega _{1}^{spin})\rightarrow
H_{2}(F/F_{k};\Omega _{1}^{spin})$ is also trivial, and the
conclusion follows.
\end{proof}

Consider the image of $\left( T_{f}^{\gamma },\phi
_{f,k+1}^{\gamma },\sigma \right) $ in $\Omega
_{3}^{spin}(F/F_{k})$ under the homomorphism $\left( \pi
_{k+1,k}\right) _{\ast }:\Omega _{3}^{spin}(F/F_{k+1})\rightarrow
\Omega _{3}^{spin}(F/F_{k}).$ This image is of course $\left(
T_{f}^{\gamma },\phi _{f,k}^{\gamma },\sigma \right) $, and since
$\left( T_{f}^{\gamma },\phi _{f,k+1}^{\gamma },\sigma \right) \in
J_{2,1}^{k+1}$, Lemma \ref{(zero)} tells us $\left( T_{f}^{\gamma
},\phi _{f,k}^{\gamma },\sigma \right) \in J_{1,2}^{k}$. By
(\ref{(Jpq)}) there is a triple $\left( M^{\prime \prime },\phi
^{\prime \prime },\sigma ^{\prime \prime }\right) \in \Omega
_{3}^{spin}((F/F_{k})^{(1)})$ such that $\left( i_{1,2}\right)
_{\ast }(M^{\prime \prime },\phi ^{\prime \prime },\sigma ^{\prime
\prime })=\left( T_{f}^{\gamma },\phi _{f,k}^{\gamma },\sigma
\right) $ in $\Omega _{3}^{spin}(F/F_{k})$ as indicated in the
following diagram.
\begin{equation*}
\begin{diagram}
                                                   &          &                         &                   &                                          &          &                         &                                                                   \\
                                                   &          &                         &                   &                                          &          & \Omega_{3}^{spin}( (F/F & _{k})^{(1)}) \ \ \ \ \ \ \ \ \ \                                  \\
                                                   &          & \Omega_{3}^{spin}( (F/F & _{k+1})           & ^{(2)}) \ \ \ \ \ \ \ \ \                &          & \dTo                    & ( M^{\prime \prime},\phi^{\prime \prime},\sigma ^{\prime \prime}) \\
                                                   &          & \dTo                    & ( M^{\prime },    & \phi^{\prime},\sigma ^{\prime }) \ \ \ \ &          &                         &                                                                   \\
                                                   &          &                         &                   & \dMapsto                                 &          &                         &                                                                   \\
\Omega_{3}^{spin}( F/F_{2k+1})                     & \rTo     & \Omega_{3}^{spin}( F/F  & _{k+1})           &                                          & \rTo     & \Omega_{3}^{spin}( F/F  & _{k}) \ \ \ \ \ \  \dMapsto                                       \\
(T_{f}^{\gamma },\phi _{f,2k+1}^{\gamma },\sigma ) & \rMapsto &                         & (T_{f}^{\gamma }, & \phi_{f,k+1}^{\gamma },\sigma )          & \rMapsto &                         & (T_{f}^{\gamma },\phi_{f,k}^{\gamma },\sigma )                    \\
\end{diagram}\bigskip
\end{equation*}

Now we use the fact that $f\in \mathcal{B}$. Recall $\omega
_{\sigma,\alpha }=\left( \pi _{\alpha }\right) _{\ast }\circ \eta
_{\sigma,2}:\mathcal{J}(2)\rightarrow \Omega _{3}^{spin}(S^{1})$
and $\left( \pi _{\alpha }\right) _{\ast }:\Omega
_{3}^{spin}(F/F_{2})\rightarrow \Omega _{3}^{spin}(S^{1})$. Since
the 1-skeleton $(F/F_{k})^{(1)}$ is homotopy equivalent to the
wedge of $2g$ circles, we have
\begin{eqnarray*}
\Omega _{3}^{spin}((F/F_{k})^{(1)}) &\cong &\Omega
_{3}^{spin}(S^{1}\vee
\cdot \cdot \cdot \vee S^{1}) \\
&\cong &\overset{2g}{\oplus }\Omega _{3}^{spin}(S^{1})
\end{eqnarray*}%
and the following commutative diagram.
\begin{equation*}
\begin{diagram}
               &                        & & \overset{2g}{\oplus }\Omega_{3}^{spin}\left(S^{1}\right) & \rTo                                  & \Omega_{3}^{spin}\left( S^{1}\right)     \\
               &                        & & \dTo_{\left( i_{1,2}\right)_{\ast }}                     &                                       & \uTo_{\left(\pi_{\alpha}\right)_{\ast }} \\
\mathcal{J}(k) & \rTo^{\eta_{\sigma,k}} & & \Omega_{3}^{spin}( F/F_{k})                              & \rTo^{\left(\pi_{k,2}\right)_{\ast }} & \Omega_{3}^{spin}( F/F_{2})              \\
\end{diagram}
\end{equation*}%
There is a basis of $H^{1}(T_{f}^{\gamma };\mathbb{Z})$ such that
for each basis element $\alpha _{i}$ the range of the homomorphism
$\omega _{\sigma,\alpha _{i}}$ corresponds to a summand of $\Omega
_{3}^{spin}((F/F_{k})^{(1)})\cong \overset{2g}{\oplus }\Omega
_{3}^{spin}(S^{1})$. Since $f\in \mathcal{B}$, $ \left( M^{\prime
\prime },\phi ^{\prime \prime },\sigma ^{\prime \prime }\right)
\in \Omega _{3}^{spin}((F/F_{k})^{(1)})$ must be trivial in each
summand of $\overset{2g}{\oplus }\Omega _{3}^{spin}(S^{1})$, and
thus it must be trivial in $\Omega _{3}^{spin}((F/F_{k})^{(1)})$.
Therefore $0=\left( i_{1,2}\right) _{\ast }(M^{\prime \prime
},\phi ^{\prime \prime },\sigma ^{\prime \prime })=\left(
T_{f}^{\gamma },\phi _{f,k}^{\gamma },\sigma \right) $ in $\Omega
_{3}^{spin}(F/F_{k})$.
\end{proof}

\end{document}